\documentclass[12pt]{article}
\usepackage{amsmath, amssymb}
\usepackage{graphicx}
\newtheorem{theorem}{Theorem}

\def\be{\begin{equation}}
\def\ee{\end{equation}}

\newcommand{\dX}{\delta X}
\newcommand{\dY}{\delta Y}
\newcommand{\dW}{\delta W}
\newcommand{\dZ}{\delta Z}
\newcommand{\e}{\textbf{e}}
\newcommand{\et}{\tilde{\textbf{e}}}

\newcommand{\Th}{\Theta}
\newcommand{\PT}{\mathcal{PT}}
\newcommand{\wt}{\tilde{w}}
\newcommand{\as}{$\alpha$-surface}
\newcommand{\bs}{$\beta$-surface}
\newcommand{\OO}{\mathcal{O}}
\newcommand{\CP}{\mathbb{CP}^1}

\newcommand{\ep}{\epsilon}

\newcommand{\koniec}{\begin{flushright}  $\Box $ \end{flushright}}

\newcommand{\R}{\mathbb{R}}
\newcommand{\hook}{{\setlength{\unitlength}{11pt}   
                   \begin{picture}(.833,.8)
                   \put(.15,.08){\line(1,0){.35}}
                   \put(.5,.08){\line(0,1){.5}}
                   \end{picture}}}
\newcommand{\C}{\mathbb{C}}

\newcommand{\Z}{\mathbb{Z}}
\newcommand{\bK}{\mathbb{K}}
\newcommand{\p}{\partial}
\newcommand{\tl}{\tilde{\lambda}}
\newcommand{\lam}{\lambda}

\def\Om{\Omega}
\newcommand{\Sm}{\mathbf{\Sigma}}

\newcommand{\M}{\mathcal{M}}

\newcommand{\W}{\mathcal{W}}
\newcommand{\tm}{\tilde{\mu}}
\newcommand{\ttau}{\tilde{\tau}}

\begin{document}
\pagestyle{plain}
\title{\vskip -70pt
\begin{flushright}
{\normalsize DAMTP-2006-90} \\
\end{flushright}
\vskip 80pt
{\bf Anti-Self-Dual Conformal Structures in Neutral Signature} \vskip 20pt}

\author{Maciej Dunajski\thanks{email m.dunajski@damtp.cam.ac.uk} \\[10pt]
and \\[10pt]
Simon West\thanks{email S.D.C.West@damtp.cam.ac.uk} \\[15pt]
{\sl Department of Applied Mathematics and Theoretical Physics} \\[5pt]
{\sl University of Cambridge} \\[5pt]
{\sl Wilberforce Road, Cambridge CB3 0WA, UK} \\[10pt]}
\date{} 
\maketitle
\begin{abstract}
We review the subject of four dimensional anti-self-dual 
conformal structures with signature $(++--)$. Both local and global 
questions are discussed. Most of the material is well known in the literature
and we present it in a way which underlines the connection with integrable 
systems.
Some of the results - e.g.
the Lax pair characterisation of the scalar--flat
K\"ahler condition and a twistor construction of a conformal 
structure with twisting null Killing vector - are new.
\end{abstract}
\maketitle
\tableofcontents
\section{Introduction}
We begin with some well-known facts from Riemannian geometry. Given
an oriented Riemannian 4-manifold $(\M,g)$, the
Hodge-$\ast$ operator is an
involution on 2-forms. This induces a decomposition
\be \label{splitting}
\Lambda^{2} = \Lambda_{+}^{2} \oplus \Lambda_{-}^{2}
\ee
of 2-forms
into self-dual
and anti-self-dual\index{anti-self-dual} components, which
only depends on the conformal class $[g]$. Now choose $g \in [g]$.
The Riemann tensor
has the index symmetry $R_{abcd}=R_{[ab][cd]}$ so can be thought of
as a map $\mathcal{R}: \Lambda^{2} \rightarrow \Lambda^{2}$.
This map decomposes under (\ref{splitting}) as follows:
\be \label{decomp}
{\mathcal R}=
\left(
\mbox{
\begin{tabular}{c|c}
&\\
$C_++\frac{s}{12}$&$\phi$\\ &\\
\cline{1-2}&\\
$\phi$ & $C_-+\frac{s}{12}$\\&\\
\end{tabular}
} \right) .
\ee
The $C_{\pm}$ terms are the self-dual and anti-self-dual parts
of the Weyl tensor, the $\phi$ terms are the
tracefree Ricci curvature, and $s$ is the scalar curvature which acts
by scalar multiplication. The Weyl tensor\index{Weyl tensor} is  conformally
invariant, so can be
thought of as being defined by the conformal structure $[g]$.
An \emph{anti-self-dual
conformal
structure} is one with $C_{+}=0$.
Such structures have
a  global twistor correspondence \cite{AHS} which has been
studied intensively; they have also been studied from a
purely analytic point of view using elliptic techniques \cite{taubes}.

What happens in other signatures? In Lorentzian signature $(+++-)$,
the Hodge-$\ast$ is not
an involution (it squares to $-\textbf{1}$ instead of $\textbf{1}$)
and there is no decomposition of 2-forms. In neutral
$(++--)$ signature,
the Hodge-$\ast$ \emph{is} an involution, and there is a
decomposition exactly as in the Riemannian
case, depending on $[g]$. Thus anti-self-dual
conformal structures exist in neutral signature. This article is
devoted to their properties.

At the level of PDEs, the difference between neutral
and Riemannian is that in the neutral case the gauge-fixed
anti-self-duality equations
are \emph{ultrahyperbolic}, whereas in the Riemannian case they are elliptic.
This results in profound differences, both locally
and globally. Roughly speaking, the neutral case is far less rigid than
the Riemannian case.
For instance, any Riemannian anti-self-dual conformal
structure must be analytic by the twistor construction. In the
neutral case there is no
general twistor construction, and in fact neutral
conformal structures are not necessarily analytic.
This lack of analyticity provides scope for rich local behaviour,
as wave like solutions exists.

Assuming symmetries in the form of Killing vectors, one
often finds that the equations reduce to \emph{integrable
systems}. Different
integrable systems can be obtained by combining symmetries with
geometric conditions for a metric in a conformal class.  The
story here in some sense parallels the case of the self-dual
Yang-Mills equations in neutral signature, where imposing
symmetries leads to many well-known integrable
systems \cite{masonwoodhouse,D_book}.

The subject of this review is the interplay between the
ultrahyperbolic differential equations, and the anti-self-duality
condition. We shall make a historical digression, and note
that both concepts arouse
separately in mid 1930s.

Indeed, the ultrahyperbolic wave equation appears naturally
in integral geometry, where  the X-ray transform\index{X-ray transform} introduced in 1938 by
John \cite{john} can be used to construct all its smooth solutions.
This takes a smooth function on $\mathbb{RP}^3$
(a compactification of $\R^3$)
and integrates it  over an oriented geodesic.
The resulting function is defined on the Grassmannian
$\mbox{Gr}_2(\R^4)$ of two-planes in $\R^4$ and satisfies the
wave equation for a flat metric in $(++--)$ signature. To see it
explicitly  consider a smooth function
$f:\R^3\longrightarrow \R$  which satisfies suitable decay
conditions at infinity. For any oriented
line $L\subset \R^3$ define
$
\psi(L)=\int_L f$, or
\be
\label{radon_john}
\psi(x,y,w,z)=\int^{\infty}_{-\infty}f(xs+z, ys-w, s)ds,
\ee
where we have chosen an explicit parametrisation of all lines which
are not perpendicular to the $x_3$ axis. The dimension of the space
of oriented lines is $4$. This is greater than the dimension of
$\R^3$, and one does not expect $\psi$ to be
arbitrary. Differentiating under the integral sign shows that
$\psi$ must satisfy  the wave equation in neutral signature
\be \label{uhwaveequation}
\frac{\partial^2\psi}{\partial x\partial w}
+\frac{\partial^2\psi}{\partial y\partial z}=0.
\ee
John has demonstrated that equation (\ref{uhwaveequation}) is the only
condition constraining the range of the integral transform in this
case, and that all smooth solutions to (\ref{uhwaveequation}) arise by
(\ref{radon_john})
from some $f$. One can regard the $X$--ray transform as the predecessor
of twistor theory. In this context $\mathbb{RP}^3$ should be regarded as a
totally real submanifold of a twistor space $\mathbb{CP}^3$. In fact
Woodhouse \cite{Wo92} showed that any local solution of
(\ref{uhwaveequation}) can be generated from a function on the real
twistor space of $\R^{2,2}$. The twistor space is the set of totally
null self-dual  2-planes and is three dimensional, so we are again
dealing with a function of three variables. To obtain the value $\psi$ at a
point $p$, one integrates $f$ over all the planes through $p$. This
was motivated by the Penrose transform with neutral reality conditions.

It is less well known that the ASD equation on Riemann curvature dates
back to the same period  as the work of John
(at least 40  years before the seminal work
of Penrose \cite{penrose} and Atyiah--Hitchin--Singer \cite{AHS}).
It arose in
the context of Wave Geometry -- a subject developed in Hiroshima
during the 1930s. Wave Geometry postulates the
existence of a privileged spinor field which in the modern
super--symmetric context  would be
called a Killing spinor. The integrability conditions come down
to the ASD condition on a Riemannian curvature of the
underlying complex space time. This condition implies vacuum Einstein
equations.
The Institute at Hiroshima where Wave Geometry had been  developed was
completely destroyed by the atomic bomb in 1945.
Two  of the survivors
wrote up the results in a book \cite{MTbook}.
In particular in \cite{old_pleban} it was shown that local
coordinates can be found such that the metric takes  a form
\be
\label{heavenly_met}
g=\frac{\p^2 \Omega}{\p x\p w}d x d w+\frac{\p^2 \Omega}{\p y\p z}d y d z
+\frac{\p^2 \Omega}{\p y\p w}d y d w+\frac{\p^2 \Omega}{\p x\p z}d x d z
\ee
and ASD vacuum condition reduces to a single PDE for one function $\Omega$
\be
\label{heavy1}
\frac{\p^2 \Omega}{\p x\p w}\frac{\p^2 \Omega}{\p y\p z}-
\frac{\p^2 \Omega}{\p x\p z}\frac{\p^2 \Omega}{\p y\p w}
=1.
\ee
This is nowadays known as the first heavenly equation after Plebanski who
rediscovered it in 1975 \cite{Pl75}.
If $(\Omega, x, y, w, z)$ are all real, the resulting
metric has neutral signature. The flat metric corresponds to $\Omega=wx+zy$.
Setting \[
\Omega=wx+zy+\psi(x,y,w,z)\]
we see that up to the linear terms in
$\psi$ the heavenly equation reduces to the ultrahyperbolic wave equation
(\ref{uhwaveequation}). Later  we shall see that the twistor method of
solving (\ref{heavy1}) is a non-linear version of John's X-Ray transform.
This concludes our historical digression.

The article is structured as follows. In Section \ref{preliminaries}
we introduce the local theory of
neutral anti-self-dual conformal structures.
It is convenient to use
spinors, which for us will be a local tool to make the  geometric
structures more transparent.
In Section \ref{integrable} we explain how neutral ASD conformal
structures are
related to Lax pairs and hence integrable systems. We review various
curvature restrictions on a metric in a conformal class (Ricci-flat,
scalar flat K\"ahler etc), and show how these can be characterised in
terms of their Lax pair.
Section \ref{symmetries} is devoted to symmetries; in this section we
make contact with many well known integrable systems.
We discuss
twistor theory in Section \ref{twistor},
explaining the differences between
the Riemannian and neutral case, and describing various twistor methods of
generating neutral ASD conformal structures.
Despite the ultrahyperbolic nature of the equations, some strong
global results have been obtained in recent years using a variety of
techniques. We discuss these in Section \ref{global}.

The subject of neutral anti-self-dual
conformal structures is rather diverse. We hope to present a
coherent overview, but the different strands will not all be woven
together. Despite this, we hope the article serves a useful purpose as a
path through the literature.

\newpage
\section{Local geometry in neutral signature}
\label{preliminaries}
\subsection{Conformal compactification}
\label{conf_comp}
We shall
start off by describing a conformal 
compactification\index{conformal compactification} of
the flat neutral metric.
Let $\R^{2, 2}$ denote  $\R^4$ with a flat $(+ + - -)$ metric. Its
natural compactification is a projective quadric in $\mathbb{RP}^5$. To
describe it explicitly consider $[{\bf x}, {\bf y}]$ as homogeneous
coordinates on $\mathbb{RP}^5$, and set ${\cal Q}=|{\bf x}|^2-|{\bf y}|^2$.
Here $({\bf x}, {\bf y})$ are vectors on $\R^3$ with its natural inner
product. The cone ${\cal Q}=0$ is projectively invariant, and
the freedom $({\bf x}, {\bf y})\sim (c{\bf x}, c{\bf y})$,
where $c\neq 0$ is fixed to set $|{\bf x}|=|{\bf y}|=1$ which is $S^2\times
S^2$. We need
to quotient this by the antipodal map $({\bf x}, {\bf y})\rightarrow
(-{\bf x}, -{\bf y})$ to obtain the conformal compactification\footnote{This
compactification can be identified with the Grassmannian $\mbox{Gr}_2(\R^4)$
arising in the John transform (\ref{radon_john}).}
\[
\overline{\R^{2, 2}}={(S^2\times S^2)}/{\Z_2}.
\]
Parametrising the double cover of this compactification by
stereographic coordinates we find that the flat metric
$|d {\bf x}|^2-|d {\bf y}|^2$ on $\R^{3, 3}$ yields
the metric
\be
\label{standard}
g_0 = 4 \frac{d\zeta d \bar{\zeta}}{(1+\zeta \bar{\zeta})^{2}} -
4 \frac{d\chi d \bar{\chi}}{(1+\chi \bar{\chi})^{2}}
\ee
on $S^2\times S^2$. To obtain the flat metric on $\R^{2. 2}$ we would instead
consider the intersection of the zero locus of ${\cal Q}$ in $\R^{3, 3,}$
with a null hypersurface $x_0-y_0=1$.

The metric $g_0$ is conformally flat and scalar flat, as the scalar
curvature is the difference between curvatures on both factors.
It is also K\"ahler with respect to the natural complex
structures on $\CP\times\CP$ with holomorphic coordinates $(\zeta,
\chi)$. In Section \ref{Tod_sub} we shall see that $g_0$ admits
nontrivial scalar--flat K\"ahler deformations \cite{tod} globally defined
on $S^2\times S^2$.

\subsection{Spinors} \label{spinors}
It is often convenient in four dimensions to use spinors\index{spinors}, 
and the
neutral
signature case is no exception.  The relevant Lie
group isomorphism in neutral signature is
\be \label{iso}
SO(2,2) \cong SL(2,\mathbb{R}) \times SL(2,\mathbb{R}) /
\mathbb{Z}_{2}.
\ee

We shall assume that the neutral four manifold
$(\M, g)$ has a spin structure.
Therefore there exist real two-dimensional vector bundles
$S, S'$  (spin-bundles) over $\M$
equipped with parallel symplectic
structures
$\epsilon, \epsilon'$ such that
$T\M\cong {S}\otimes {S'}$
is a  canonical bundle isomorphism, and
\[
g(v_1\otimes w_1,v_2\otimes w_2)
=\epsilon(v_1,v_2)\epsilon'(w_1,w_2)
\]
for $v_1, v_2\in \Gamma(S)$ and $w_1, w_2\in \Gamma(S')$.
The two-component spinor notation \cite{PR86}  will used in
the paper.
The spin bundles $S$ and $S'$
inherit connections from the Levi-Civita
connection such that $\ep$, $\epsilon'$ are covariant constant. We
use the standard convention in which spinor indices are capital
letters, unprimed for sections of $S$ and primed for sections of $S'$.
For example $\mu_{A}$ denotes a section of $S^{*}$, the dual of $S$,
and $\nu^{A'}$ a
section of $S'$.

The symplectic structures on spin spaces $\epsilon_{AB}$
and $\epsilon_{A'B'}$ (such that
$\epsilon_{01}=\epsilon_{0'1'}=1$) are used to raise and
lower indices. For
example given a section $\mu^{A}$ of $S$ we define a section of $S^{*}$
by $\mu_{A}:=\mu^{B}\epsilon_{BA}$.

Spin dyads $(o^A, \iota^{A})$ and $(o^{A'}, \iota^{A'})$ span $S$
and $S'$ respectively.
We denote a  normalised null tetrad\index{null tetrad} of vector fields on $\M$ by
\[
\e_{AA'}=
\left (
\begin{array}{cc}
\e_{00'}&\e_{01'}\\
\e_{10'}&\e_{11'}
\end{array}
\right ).
\]
This tetrad is determined by the choice of spin dyads in the sense that
\[
o^{A}o^{A'}\e_{AA'}=\e_{00'},\;
\iota^{A}o^{A'}\e_{AA'}=\e_{10'},\;
o^{A}\iota^{A'}\e_{AA'}=\e_{01'},\;
\iota^{A}\iota^{A'}\e_{AA'}=\e_{11'}.
\]
The dual tetrad of one-forms by $\e^{AA'}$ determine the metric by
\be
\label{themetric}
g=
\epsilon_{AB}\epsilon_{A'B'}\e^{AA'}\otimes \e^{BB'}
=2(\e^{00'}\odot \e^{11'}-\e^{10'}\odot \e^{01'})
\ee
where $\odot$ is the symmetric tensor product.
With indices, the above formula\footnote{Note that we
drop the prime on $\epsilon'$ when using indices, since it is
already distinguished from $\epsilon$ by the primed indices.}
for $g$ becomes
$g_{ab}=\epsilon_{AB}\epsilon_{A'B'}$.

The local basis $\Sm^{AB}$ and $\Sm^{A'B'}$ of  spaces of ASD and SD
two-forms are defined by
\be\label{sdforms}
\e^{AA'}\wedge \e^{BB'}=\epsilon^{AB}\Sm^{A'B'}+\epsilon^{A'B'}\Sm^{AB}.
\ee

A vector $V$
be decomposed as $V^{AA'}\e_{AA'}$,
where $V^{AA'}$ are the components of $V$ in the basis. Its norm
is given by $\text{det} (V^{AA'})$, which is unchanged under
multiplication of the matrix $V^{AA'}$ by elements of $SL(2,\R)$ on
the left and right
\[
V^{AA'}\longrightarrow \Lambda^A_BV^{BB'}\Lambda^{A'}_{B'}, \qquad
\Lambda\in SL(2, \R),\quad \Lambda^{'}\in SL(2, \R)^{'}
\]
giving (\ref{iso}). The quotient by
$\mathbb{Z}_{2}$ comes from the fact that multiplication on the left
and
right by  $-\mathbf{1}$ leaves $V^{AA'}$ unchanged.

Spinor notation is particularly useful for describing null
structures. A vector $V$ is null when $\text{det} (V^{AA'})=0$, so
$V^{AA'}=\mu^{A}\nu^{A'}$ by linear algebra. In invariant language,
this says that a vector
$V$ is null iff $V=\mu \otimes \nu$ where $\mu, \nu$ are sections of
$S, S'$.

The decomposition of a 2-form into self-dual and anti-self-dual
parts is straightforward in spinor notation. Let $F_{AA'BB'}$ be a
2-form in indices.
Now
\begin{eqnarray*}
F_{AA'BB'} &=&
F_{(AB)(A'B')}+F_{[AB][A'B']}+F_{(AB)[A'B']}+F_{[AB](A'B')}\\
&=& F_{(AB)(A'B')} + c \epsilon_{AB} \ep_{A'B'} +
\phi_{AB}\ep_{A'B'} + \psi_{A'B'}\ep_{AB}.
\end{eqnarray*}
Here we have used the fact that in two dimensions there is a
unique anti-symmetric matrix up to scale, so whenever an
anti-symmetrized pair of spinor indices occurs we can substitute
a multiple of $\epsilon_{AB}$ or $\ep_{A'B'}$ in their place. Now
observe that the first two
terms are incompatible with $F$ being a 2-form, i.e.
$F_{AA'BB'}=-F_{BB'AA'}$. So we obtain
\be \label{2form}
F_{AA'BB'} = \phi_{AB}\ep_{A'B'} + \psi_{A'B'}\ep_{AB},
\ee
where $\phi_{AB}$ and $\psi_{A'B'}$ are symmetric. This is precisely
the decomposition of $F$ into self-dual and anti-self dual parts.
Which is which depends on the choice of volume form; we choose
$\psi_{A'B'}\ep_{AB}$ to be the self-dual part. Invariantly, we have
\be
\Lambda^{2}_{+} \cong S'^{*} \odot S'^{*}, \ \ \ \ \
\Lambda^{2}_{-} \cong S^{*} \odot S^{*}.
\ee

\subsection{$\alpha$ and $\beta$ planes} \label{alphabeta}
Suppose at a point $x\in {\cal M}$
 we are given a spinor $\nu^{A'} \in S'_{x}$.
A two-plane $\Pi_{x}$ is defined by all vectors of the form
$V^{AA'}=\mu^{A}\nu^{A'}$, with varying $\mu^{A}\in S$. Now suppose $V,W \in
\Pi_{x}$. Then
$g(V,W)=\nu^{A'}\nu^{B'}\ep_{A'B'}\kappa^{A}\mu^{B}\epsilon_{AB}=0$
since $\ep_{A'B'}$ is antisymmetric. Therefore we say the two-plane is
\emph{totally
null}. Furthermore, the 2-form $V_{[a}W_{b]}$ is proportional to
$\mu_{A'}\mu_{B'}\epsilon_{AB}$; i.e. the two-plane is self-dual. In
summary, a spinor in $S$ defines a totally null self-dual
two-plane, which is called an \emph{$\alpha$-plane}\index{$\alpha$-plane}. 
Similarly a
spinor in $S$ defines a totally null anti-self-dual two-plane,
called a \emph{$\beta$-plane}\index{$\beta$-plane}.

\subsection{Anti-self-dual conformal structures in spinors}

A neutral conformal structure $[g]$ is an equivalence
class of neutral signature metrics, with the equivalence relation $g
\sim e^{f} g$ for any function $f$. Another way of viewing such a structure is
as a line-bundle valued neutral metric; we will not need this
description because in most cases we will be working with
particular metrics within a conformal class.

Choose a $g \in [g]$. Then there is a Riemann tensor, which possesses
certain symmetries under permutation of indices. In the same way that
we deduced (\ref{2form}) for the decomposition of a 2-form in
spinors, the Riemann tensor decomposes as \cite{PR86}
\begin{eqnarray*}
R_{AA'BB'CC'DD'} &=& C_{ABCD}\ep_{A'B'}\ep_{C'D'} +
C_{A'B'C'D'}\ep_{AB}\ep_{CD} \\ && + \phi_{ABC'D'}\ep_{A'B'}\ep_{CD}
+ \phi_{A'B'CD}\ep_{AB}\ep_{C'D'}  \\ && + \frac{s}{12}
(\ep_{AC}\ep_{BD}\ep_{A'B'}\ep_{C'D'} +
\ep_{AB}\ep_{CD}\ep_{A'D'}\ep_{B'C'}).
\end{eqnarray*}
This is  the spinor version of (\ref{decomp}). Here
$C_{A'B'C'D'}$, $C_{ABCD}$  are totally symmetric, and correspond to
$C_{+}$, $C_{-}$ in
(\ref{decomp}).
The spinor $\phi_{A'B'CD}$ is symmetric in its pairs of indices, and
corresponds to $\phi$ in (\ref{decomp}).
 An
anti-self-dual conformal structure
is one for which $C_{A'B'C'D'}=0$. In the next section we explain the
geometric significance of this condition in more detail.
It is appropriate here to recall the Petrov-Penrose classification
\index{Petrov-Penrose classification} 
\cite{PR86} of the algebraic type of a Weyl tensor. In split
signature this applies separately to $C_{ABCD}$ and
$C_{A'B'C'D'}$. In our case $C_{A'B'C'D'}=0$ and we are concerned with the
algebraic type of $C_{ABCD}$. One can form a
real polynomial of fourth
order $P(x)$ by defining $\mu^A = (1,x)$ and setting $P(x) = \mu^A
\mu^B \mu^C \mu^D C_{ABCD}$. The Petrov-Penrose classification refers
to the position of roots of this polynomial, for example if there are
four repeated roots then we say $C_{ABCD}$ is type N\index{type N}. 
If there is a
repeated root the metric is called \emph{algebraically special}\
index{algebraically special}.
There are additional complications in the split signature case
\cite{Law} arising from the fact that real polynomials may not have
real roots.

\section{Integrable systems and Lax pairs} \label{integrable}
In this section we show how anti-self-dual conformal structures are
related to integrable
systems and Lax pairs\index{Lax pair}. Let $g\in [g]$ and
let $\nabla$ denote the Levi--Civita
connection on $\M$. This connection induces spin 
connections\index{spin connection} on
spin bundles which we also denote $\nabla$. Let us consider $S'$.
The connection coefficients $\Gamma_{AA'B'}^{\
\ \ \ \ \ C'}$ of $\nabla$ are defined by
$$
\nabla_{AA'} \mu^{C'} = \e_{AA'} (\mu^{C'}) + \Gamma_{AA'B'}^{\ \ \ \
\ \ C'}¥ \mu^{B'},
$$
where $\mu^{A'}$ is a section of $S'$ in coordinates determined by
the basis $\e_{AA'}$. The $\Gamma_{AA'B'}^{\ \ \ \
\ \ C'}$ symbols can  be calculated in terms of the Levi-Civita
connection symbols. They can also be read off directly from
the Cartan equations $d \e^{AA'}=\e^{BA'}\wedge {\Gamma_B}^A+
\e^{AB'}\wedge {\Gamma_{B'}}^{A'}$, where ${\Gamma_{B'}}^{C'}=
\Gamma_{AA'B'}^{\ \ \ \
\ \ C'}{\bf e}^{AA'}$.
See \cite{PR86} for details. Now given a connection
on a vector bundle, one can lift a vector field on the base to
a horizontal vector field on the total space. We follow standard
notation and denote the local coordinates of $S'$ by $\pi^{A'}$. Then
the horizontal lifts $\et_{AA'}$ of $\e_{AA'}$ are given explicitly by
$$
\et_{AA'} := \e_{AA'} + \Gamma_{AA'B'}^{\ \ \ \
\ \ C'} \pi^{B'}\frac{\p}{\p \pi^{C'}}.
$$
Now we can state a seminal result of Penrose:
\begin{theorem}{\cite{penrose}} \label{penrose}
Given a neutral metric $g$, define a two dimensional distribution on
$S'$ by $\mathcal{D} = \text{span} \{
L_{0}, L_{1} \}$, where
\be
\label{laxpair}
L_{A}:=\pi^{A'} \et_{AA'}.
\ee
Then $\mathcal{D}$ is integrable iff $g$ is anti-self-dual.
\end{theorem}
So when $g$ is ASD, $S'$ is foliated by surfaces. Since the $L_{A}$
are homogeneous in the
$\pi^{A'}$ coordinates, $\mathcal{D}$ defines a distribution
on $PS'$, the projective version of $S'$.

The push down of
$\mathcal{D}$ from a point $\pi^{A'}=\nu^{A'}$ in a fibre of $S'$ to the
base is the $\alpha$-plane defined by $\nu^{A'}$, as explained in
Section \ref{alphabeta}.
So the content of Theorem \ref{penrose}  is that
$g$ is ASD iff any $\alpha$-plane is tangent to an 
$\alpha$-surface\index{$\alpha$-surface},
i.e. a surface that is totally null and self-dual at every point.
Any such $\alpha$-surface lifts to a unique surface in $PS'$,
or a one parameter
family of surfaces in $S'$.
\subsection{Curvature restrictions and their Lax pairs}
\label{specialmetrics}
A more recent interpretation of Theorem \ref{penrose} is to regard $L_{A}$
as a Lax pair for the ASD conformal structure. Working on $PS'$, with
inhomogeneous fibre coordinate $\lambda=\pi^{1'}/\pi^{0'}$, the condition that
$\mathcal{D}$ commutes is the compatibility condition for the pair of
linear equations
\begin{eqnarray*}
L_{0}f=(\et_{00'} + \lambda \et_{01'}) f &=& 0 \\
L_{1}f=(\et_{10'} + \lambda \et_{11'}) f &=& 0
\end{eqnarray*}
to have a solution $f$  for all $\lambda \in \R$, where $f$ is a function on
$PS'$. In integrable systems language, $\lambda$ is the spectral
parameter\index{spectral parameter}.

Here we describe various conditions that one can place on a metric $g
\in [g]$ on top of anti--self--duality. This provides a more direct link
with integrable systems as in each case described below one can choose
a spin frame, and local coordinates to reduce the special ASD condition
to an integrable scalar PDE with corresponding Lax pair.


\subsubsection{Pseudo-hyperhermitian structures}
This is the neutral analogue of Riemannian hyperhermitian geometry.
The significant point for us is that in four dimensions,
pseudo-hyperhermitian\index{pseudo-hyperhermitian metrics} 
metrics (defined below) are necessarily
anti-self-dual.

Consider
a structure $(\M,I,S,T)$, where $\M$ is a 4-dimensional manifold and $I,S,T$
are anti-commuting endomorphisms of the tangent bundle satisfying
\be  \label{paraquat}
S^{2}=T^{2}=\textbf{1}, \ \ I^{2}=-\textbf{1}, \ \ ST=-TS=\textbf{1}.
\ee
This is called the algebra of para-quaternions \cite{ivanov} or split
quaternions \cite{DJS}. Consider the
hyperboloid of almost complex structures on $\M$ given by
$aI+bS+cT$, for $(a,b,c)$ satisfying $a^{2}-b^{2}-c^{2}=1$.
If each of these almost complex
structures is integrable, 
we call $(\M,I,S,T)$ a \emph{pseudo-hypercomplex}
manifold.

So far we have not introduced a metric. A natural
restriction on a metric given a pseudo-hypercomplex structure is to
require it to be hermitian with respect
to each of the complex structures. This is equivalent to the
requirement:
\be \label{hermitian}
g(X,Y)=g(IX,IY)=-g(SX,SY)=-g(TX,TY),
\ee
for all vectors $X,Y$. A metric satisying (\ref{hermitian}) must be neutral. To see this
consider the endomorphism $S$, which squares to the identity. Its
eigenspaces decompose into $+1$ and $-1$ parts. Any eigenvector must
be null from (\ref{hermitian}). So choosing an eigenbasis one can
find 4 null vectors, from which it follows that the metric is
neutral. Given a pseudo-hypercomplex manifold, we call a metric satisfying
(\ref{hermitian}) a pseudo-hyperhermitian metric. There exists a unique 
conformal structure  associated to a pseudo-hypercomplex in this way 
\cite{gueo}.


As mentioned above, it turns out that pseudo-hyperhermitian metrics
are necessarily anti-self-dual. One way to formulate this is via the
Lax pair formalism as follows:
\begin{theorem}{\cite{duna:twisted}} \label{twistedphoton}
    Let $\e_{AA'}$ be four independent
    vector fields on a
    four dimensional real manifold $\M$. Put
    $$
    L_{0}=\e_{00'}+\lambda \e_{01'}, \ \ L_{1}=\e_{10'}+\lambda
    \e_{11'}.
    $$
    If
    \be \label{hhcommutator}
    [L_{0},L_{1}]=0
    \ee
    for every value of a parameter $\lambda$, then  $g$ given by
(\ref{themetric})
a pseudo-hyperhermitian metric on $\M$. Given any
    four-dimensional pseudo-hyperhermitian metric there
exists a null tetrad such that (\ref{hhcommutator}) holds.
 \end{theorem}
Interpreting
$\lambda$ as the projective primed spin coordinate as in Section
\ref{integrable}, we see that a pseudo-hyperhermitian metric must be ASD
from Theorem \ref{penrose}.
Theorem \ref{twistedphoton} characterises
pseudo-hyperhermitian metrics as those which possess a Lax pair
containing no $\p_{\lambda}$ terms.



 We shall now discuss the local formulation of the
pseudo-hyperhermitian condition as a PDE. Expanding
equation (\ref{hhcommutator}) in powers of $\lambda$ gives
\be \label{hhequations}
[\e_{A0'},\e_{B0'}]=0, \ \ \
[\e_{A0'},\e_{B1'}]+[\e_{A1'},\e_{B0'}]=0, \ \ \ [\e_{A1'},\e_{B1'}]=0.
\ee
It follows from (\ref{hhequations}), using the
Frobenius theorem and the Poincar\'e lemma, that one can
choose coordinates $(p^{A},w^{A})$ , ($A=0,1$), in which $\e_{AA'}$
take the form
\[
\e_{A0'}=\frac{\p}{\p p^{A}}, \qquad \e_{A1'}=\frac{\p}{\p w^{A}} -
\frac{\p \Theta^{B}}{\p p^{A}} \frac{\p}{\p p^{B}},
\]
where $\Theta^{B}$ are a pair of functions satisfying a system
of coupled non-linear
ultra-hyperbolic PDEs.
\be  \label{hhheavenly}
\frac{\p^{2} \Theta_{C}}{\p p_{A} \p w^{A}}+\frac{\p\Theta_{B}}{\p
p^{A}}\frac{\p^{2} \Theta_{C}}{\p p_{A} \p p_{B}}=0.
\ee
Note the indices here are not spinor indices, they are simply a
convenient way of labelling coordinates and the functions
$\Theta^{A}$. We raise and lower them using the standard antisymmetric matrix
$\ep_{AB}$, for example $p_{A}:=p^{B}\ep_{BA}$, and the summation
convention is used.

\subsubsection{Scalar--flat K\"ahler structures}

Let $(\M, g)$ be an ASD four manifold and let $J$ be a (pseudo) complex
structure such that the corresponding fundamental two--form is closed.
This ASD K\"ahler condition implies that $g$ is scalar flat, and conversely
all scalar flat K\"ahler four manifolds are ASD  \cite{Derdzinski}.

In this subsection we shall show that in the 
scalar--flat K\"ahler\index{scalar--flat K\"ahler metrics} case the
spin frames can be chosen  so that the Lax pair
(\ref{laxpair}) consists of
volume-preserving vector fields on $\M$ together with two functions
on $\M$.
The following theorem has been obtained in a joint work of Maciej Przanowski
and the first author. We shall formulate and prove it in the
holomorphic category which will allow both neutral and Riemannian
real slices.
\begin{theorem}
\label{d01}
Let
${\e}_{AA'}=({\e}_{00'}, {\e}_{01'},
{\e}_{10'}, {\e}_{11'})$
be four independent holomorphic vector
fields on a four-dimensional complex manifold $\M$ and let
$f_1, f_2: \M \rightarrow \C$ be two holomorphic function.
Finally, let $\nu$
be a nonzero holomorphic four-form. Put
\be
\label{sde}
L_0={\e}_{00'}+\lam{\e}_{01'} -f_0\lam^2\frac{\p}{\p \lam},
\qquad L_1={\e}_{10'}+\lam{\e}_{11'}-
f_1\lam^2\frac{\p}{\p \lam}.
\ee
Suppose that for every $\lambda \in \mathbb{CP}^1$
\be
\label{podstawka}
[L_0, L_1]=0,\qquad{\cal L}_{L_A}{\nu}=0,
\ee
where ${\cal L}_V$ denotes the Lie derivative.
Then \[
\widehat\e_{AA'} = c^{-1}{\e}_{AA'},\qquad
\mbox{where}\qquad
c^2:={\nu}({\e}_{00'},  {\e}_{01'},
{\e}_{10'}, {\e}_{11'}),
\]
is a null-tetrad for an ASD K\"ahler metric. Every such metric locally
arises in this way.
\end{theorem}
{\bf Proof.}
First assume that there exists a tetrad $\e_{AA'}$ and two
functions $f_A=(f_0, f_1)$ such that equations (\ref{podstawka}) are satisfied.
For convenience write down equations $[L_0, L_1]=0$ in full
\be
\label{eq1}
[\e_{A0'}, \e_{B0'}]=0,
\ee
\be
\label{eq2}
[\e_{A0'}, \e_{B1'}]+[\e_{A1'}, \e_{B0'}]=0,
\ee
\be
\label{eq3}
[\e_{A1'}, \e_{B1'}]=\epsilon_{AB}f^C\e_{C1'},
\ee
\be
{\e^A}_{0'}f_A=0,
\ee
\be
{\e^A}_{1'}f_A=0.
\ee
Define the almost complex structure $J$ by
\[
J(\e_{A1'})=-i\e_{A1'}, \qquad
J(\e_{A0'})=i\e_{A0'}.
\]
Equations (\ref{eq1}), and (\ref{eq3}) imply that 
this complex structure\index{complex structure} is
integrable. Let $g$ be a metric corresponding to  $\widehat{\e}_{AA'}$
by (\ref{themetric}).
To complete this part of the proof we need to show that
a fundamental two-form $\omega$ defined by $\omega(X, Y)=g(X, JY)$ is closed.
First observe that
\[
\omega=\frac{\p}{\p \lam}(\nu(L_0, L_1, .\;, .\;))|_{\lam=0}.
\]
It is therefore enough to prove that $\Sm=\nu(L_0, L_1, .\;., .\;)$
is closed for each fixed  $\lam$. We shall establish this fact using
equations (\ref{podstawka}), and $d \nu=0$. Let us calculate
\begin{eqnarray*}
d\Sm&=&d (\nu(L_0, L_1, .\;, .\;))=d (L_0\hook(\nu( L_1, .\;, .\;, .\;)))\\
&=&{\cal L}_{L_0}(\nu(L_1, .\;, .\;, .\;))-
L_0\hook(d\nu( L_1, .\;, .\;, .\;))\\
&=&[L_0, L_1]\hook \nu+L_1\hook {\cal L}_{L_0}(\nu)-
L_0\hook (L_1 \hook d\nu)\\
&=&-L_0\hook({\cal L}_{L_1}\nu- L_0\hook(L_1\hook  d(\nu)))=0.
\end{eqnarray*}
Therefore $\omega$ is closed which in the case of integrable $J$
also implies $\nabla\omega =0$ \cite{kob}.\\
{\bf Converse.}
The metric $g$ is K\"ahler, therefore there exist local coordinates
$(w^A,\wt^A)$ and a complex valued function $\Om=\Om(w^A, \wt^A)$
such that $g$ is given by
\be
\label{kmetric}
g= \frac{\p^2\Omega}{\p w^A\p \wt^B}dw^Ad\wt^B.
\ee
Choose a spin frame $(o_{A'}, \iota_{A'})$ such that
the tetrad of vector fields $\e_{AA'}$  is
\[
\e_{A0'}=o^{A'}\e_{AA'}=\frac{\p}{\p
w^{A}},\qquad
\e_{A1'}=\iota^{A'}\e_{AA'}=
\frac{\p^2 \Om}{\p w^A \p\wt^B}\frac{\p}{\p \wt_{B}}.
\]
The null tetrad for the metric
(\ref{kmetric}) is $\hat{\e}_{AA'}=G^{-1}\e_{AA'}$,
where
\be
\label{determinant}
G=\det(g)=\frac{1}{2}\frac{\p^2\Om}{\p w_A\p \wt_B}
\frac{\p^2\Om}{\p w^A\p \wt^B}.
\ee
The Lax pair (\ref{laxpair}) is
\[
L_A=\frac{\p}{\p w^A}-
\lam\frac{\p^2 \Om}{\p w^A \p\wt^B}\frac{\p}{\p \wt_{B}}
+l_A\frac{\p}{\p \lam}.
\]
Consider the Lie bracket
\begin{eqnarray*}
[L_0, L_1]&=&\l^2\frac{\p^2 \Om}{\p w^A \p\wt^B}
\frac{\p^3 \Om}{\p w_A \p\wt_B\p\wt^{C}}\frac{\p}{\p \wt_{C}}
+l^A\frac{\p^2 \Om}{\p w^A \p\wt^B}\frac{\p}{\p \wt_{B}}  \\
&+&
\Big(\frac{\p l^A}{\p w^A}
-\lam\frac{\p^2 \Om}{\p w^A \p\wt^B}\frac{\p l^A}{\p \wt_{B}}
+l_A\frac{\p l^A}{\p \lam}\Big)
\frac{\p}{\p \lam}.
\end{eqnarray*}
The ASD condition is equivalent to integrability of the distribution
$L_A$, therefore
\[
[L_A, L_B]=\epsilon_{AB}\alpha^CL_C
\]
for some $\alpha^C$.
The lack of $\p/\p w^A$ term
in the Lie bracket above implies $\alpha^C=0$.
Analysing other terms
we deduce the existence of
$f=f(w^A, \wt^A)\in \ker\square$ such that $l_A=\lam^2\p f/\p w^A $, and
\be
\label{nk11}
\frac{\p^2\Om}{\p w^A\p \wt^B}\frac{\p}{\p \wt^C}
\Big(\frac{\p^2\Om}{\p w_A\p \wt_B}\Big)
=\frac{\p f}{\p w_A}\frac{\p^2\Om}{\p w^A\p \wt^C}.
\ee
\koniec
The real-analytic $(++--)$ slices are obtained if
$\e_{AA'}, \nu, f_1, f_2$ are all real. In this case
we alter our definition of $J$ by
\[
J(\e_{A1'})=-\e_{A1'}, \qquad
J(\e_{A0'})=\e_{A0'}.
\]
Therefore $J^2={\bf 1}$, and $g$ is pseudo-K\"ahler.

In the Euclidean case the quadratic-form $g$ and the complex
structure
\[
J=i(\e^{A0'}\otimes\e_{A0'}-\e^{A1'}\otimes\e_{A1'})
\]
are real but the vector fields
$\e_{AA'}$ are complex.

As a corollary from the last theorem we can deduce a formulation
of the scalar--flat K\"ahler condition \cite{Park92} .
Scalar-flat-K\"ahler metric are locally given by
(\ref{kmetric})
where $\Om(w^A, \wt^A)$ is a solution to a 4th order PDE
(which we write as a  system of two second order PDEs ):
\be
\label{nk1new}
\frac{\p f}{\p w^A}=\frac{\p^2 \Omega}{\p w^A\p \wt^B}
\frac{\p \ln{G}}{\p \wt_B},
\ee
\be
\label{nk2}
\square f=\frac{\p^2\Om}{\p w^A\p \wt^B}\frac{\p^2 f}{\p w_A\p \wt_B}=0.
\ee
Moreover {(\ref{nk1new},\ref{nk2})} arise as an integrability condition
for the linear system $L_0\Psi=L_1\Psi=0$, where
$\Psi=\Psi(w^A, \wt^A, \lam)$
and
\be
\label{sflax}
L_A= \frac{\p}{\p
w^A}-\lam\frac{\p^2 \Om}{\p w^A \p\wt^B}\frac{\p}{\p \wt_{B}}
+\lam^2\frac{\p f}{\p w^A}\frac{\p}{\p \lam}.
\ee

To see this note that
in the Proof of Theorem \ref{d01} we have demonstrated that
$f\in \ker\square$. In the adopted coordinate system
\[
\square=\frac{\p^2 \Om}{\p w^A \p\wt^B}\frac{\p^2}{\p w_A\p \wt_{B}},
\]
which gives (\ref{nk1}).
Solving the algebraic system (\ref{nk11}) for $\p f/\p w ^A$
yields (\ref{nk2}).
\koniec

\subsubsection{Null-K\"ahler structures}
A null-K\"ahler structure\index{null-K\"ahler structures} on a 
real four-manifold $\M$
consists of an inner
product $g$ of signature $(++--)$ and a real rank-two endomorphism $
N:T\M \rightarrow T \M$ parallel
with respect to this inner product such that
\[
N^2=0,\qquad\mbox{and}\qquad g(NX, Y)+g(X, NY)=0
\]
for all $X, Y \in T \M$.
The isomorphism ${\Lambda^2}_+({\cal M})\cong {\mbox{Sym}}^2(S')$ between
the bundle of self-dual two-forms and the symmetric tensor product of
two spin bundles implies that the  existence of a null--K\"ahler structure is
in four dimensions equivalent to the existence of a parallel real spinor.
The Bianchi identity implies the vanishing of the curvature scalar.

In \cite{B00} and \cite{nullkahler} it was shown that null--K\"ahler structures are
locally given by one arbitrary  function of four variables, and
admit a canonical form\footnote{The local form (\ref{nkmetric}) is a
special case of Walker's  canonical form
of a neutral metric which admits a two--dimensional distribution
which is parallel and null \cite{Walker}.
Imposing more restrictions
on Walker's metric leads to examples of conformally Osserman structures,
i.e. metrics for which the eigenvalues of the operator
$Y^a\rightarrow C^a_{bcd}X^bY^cX^d$ are constant on the unit
pseudosphere $\{ X\in T{\cal M}, g(X, X)=\pm 1\}$. These metrics are
all SD or ASD according to \cite{osserman}.}
\be
\label{nkmetric}
g= d w d x+ d z  d y-\Th_{xx} d z^2-\Th_{yy} d w^2+2\Th_{xy} d w d z,
\ee
with
$N= d w\otimes\p/\p y- d z\otimes\p/\p x$.

Further conditions can be imposed on the curvature of $g$
to obtain non--linear PDEs for the
potential function $\Theta$.
Define
\be
\label{nk1}
f:=\Theta_{wx}+\Theta_{zy}+\Theta_{xx}\Theta_{yy}-\Theta_{xy}^2.
\ee
\begin{itemize}
\item The Einstein condition implies that
\[
f=xP(w, z)+yQ(w, z)+R(w, z),
\]
where $P, Q$ and $R$ are arbitrary functions of $(w, z)$.
In fact the number of the arbitrary functions
can be reduced down to one by redefinition of
$\Th$ and the coordinates.
This is the hyper--heavenly equation
of Pleba\'nski and Robinson \cite{PR76} for non--expanding metrics of type
$[N]\times$[Any]. (Recall that $({\cal M},  g)$ 
is called hyper--heavenly\index{hyper--heavenly metric} if
the self--dual Weyl spinor is algebraically special).
\item The conformal anti--self--duality (ASD) condition implies
a 4th order PDE for $\Theta$
\be
\label{md}
\square f=0,
\ee
where $\square$ is the Laplace--Beltrami operator defined by the metric $g$.
This equation is integrable:
It admits a Lax pair
\begin{eqnarray*}
L_0&=&(\p_w-\Th_{xy}\p_y+\Th_{yy} \p_x)-\lam\p_y+f_y\p_{\lam},\nonumber\\
L_1&=&(\p_z+\Th_{xx}\p_y-\Th_{xy} \p_x)+\lam\p_x+f_x\p_{\lam}.
\end{eqnarray*}
and its solutions can in principle be found
by twistor methods \cite{nullkahler}, or the dressing 
approach \cite{druma}.

\item Imposing both conformal ASD and Einstein condition implies (possibly
after a redefinition of $\Theta$) that $f=0$, which yields the
celebrated second heavenly equation of Pleba\'nski \cite{Pl75}
\be
\label{secondh}
\Th_{wx}+\Th_{zy}+\Th_{xx}\Th_{yy}-\Th_{xy}^2=0.
\ee
\end{itemize}
\[
\mbox{Null K\"ahler}\begin{array}{rcccl}
&& \mbox{ASD}  &&\\
&\nearrow&&\searrow &\\
&&&&\\
&\searrow&&\nearrow &\\
&&\mbox{Einstein}  &&
\end{array}
\mbox{Pseudo hyper--K\"ahler}.
\]

\subsubsection{Pseudo-hyperk\"ahler structures} \label{phk}
Suppose we are given a pseudo hypercomplex structure as defined in
the previous section, i.e. a two dimensional hyperboloid of integrable complex
structures. In the previous section we defined a pseudo-hyperhermitian metric to be a metric
that is hermitian with respect to each complex structure in the
family. If we further require that the 2-forms
\be \label{twoforms}
\omega_{I}(.,.)=g(.,I.), \ \ \ \omega_{S}(.,.)=g(.,S.), \ \ \
\omega_{T}(.,.)=g(.,T.),
\ee
be closed, we call say $g$ is pseudo-hyperk\"ahler
\index{pseudo-hyperk\"ahler metrics}. These define
three symplectic forms, and Hitchin has termed such structures
hypersymplectic\footnote{Other terminology includes neutral
hyperk\"ahler \cite{kam:kodaira} and hyper-parak\"ahler \cite{ivanov}.}
\cite{hitchin:HS}\index{hypersymplectic structure}.

It
follows from similar arguments to those in standard Riemannian
K\"ahler geometry that $(I,S,T)$ are covariant constant, and hence so
are $\omega_{I}, \omega_{S}, \omega_{T}$.

As in the Riemannian case, pseudo-hyperk\"ahler metrics are
equivalent to Ricci-flat anti-self-dual metrics. One can deduce this
by showing that the 2-forms (\ref{twoforms}) are self-dual, and
since they are also covariant constant there exists a basis of
covariant constant primed spinors. Then using the spinor Ricci
identities one can deduce anti-self-duality and Ricci-flatness. See
for details.

The Lax pair formulation for a pseudo-hyperk\"ahler metric is as
follows:
\begin{theorem}{\cite{ashtekar,masonnewman}} \label{hklax}
    Let $\e_{AA'}$ be four independent
    vector fields on a
    four dimensional real manifold $\M$, and $\nu$ be a 4-form. Put
    $$
    L_{0}=\e_{00'}+\lambda \e_{01'}, \ \ L_{1}=\e_{10'}+\lambda
    \e_{11'}.
    $$
    If
    \be \label{hhcommutator1}
    [L_{0},L_{1}]=0
    \ee
    for every $\lambda \in \mathbb{RP}^{1}$, and
    \be \label{volpreserving}
    \mathcal{L}_{L_{A}} \nu = 0,
    \ee
    then $f^{-1} \e_{AA'}$ is a
    null tetrad for a pseudo-hyperk\"ahler metric on $\M$, where
    $f^{2}=\nu (\e_{00'},\e_{0'},\e_{10'},\e_{11'})$. Given any
    four-dimensional pseudo-hyperk\"ahler metric such a null tetrad
    and 4-form exists.
\end{theorem}
The extra volume preserving condition (\ref{volpreserving}) distinguishes
this from Theorem \ref{twistedphoton}. Alternatively
Theorem \ref{hklax} arises as a special case of Theorem \ref{d01} with
$f_A=0$.


\vspace{6pt}
\textbf{The Heavenly Equations.} It was shown by
Pleba\'nski \cite{Pl75} that one can always put a pseudo-hyperk\"ahler
metric into the form  (\ref{heavenly_met}), where $\Omega$ satisfies
the \emph{first Heavenly equation}\index{ Heavenly equations} (\ref{heavy1}).
The function $\Omega$
can be interpreted as the K\"ahler potential\index{K\"ahler potential}
for one of the complex structures.
Pleba\'nski also gave the alternative local form. The metric is given by
(\ref{nkmetric}), and the potential $\Theta$ satisfies
the \emph{second Heavenly equation} (\ref{secondh}).
The Heavenly equations are
non-linear ultrahyperbolic equations.
These formulations are  convenient for understanding local
properties of pseudo--hyperk\"ahler metrics, as they only depend on a
single function satisfying a single PDE.

\section{Symmetries} \label{symmetries}
By a symmetry of a metric, we mean a conformal
Killing vector, i.e. a vector field
$K$ satisfying
\be \label{KV}
\mathcal{L}_{K} g = c \ g,
\ee
where $c$ is a function. If $c$ vanishes, $K$ is called a \emph{pure}
Killing vector\index{Killing vector}, otherwise it 
is called a \emph{conformal} Killing
vector\index{conformal Killing vector}. 
If $c$ is a nonzero constant $K$ is called a
\emph{homothety}\index{homothety}. If we are dealing with a 
conformal structure $[g]$, a
symmetry is a vector field $K$ satisfying (\ref{KV}) for some $g \in
[g]$. Then $(\ref{KV})$ will be satisfied for any
$g \in [g]$, where the function $c$ will depend on the
choice of $g \in [g]$. Such a $K$ is referred to as a conformal
Killing vector for the conformal structure.

In neutral signature there are two types of Killing vectors: non-null
and null. Unlike in the Lorentzian case where non-null vectors can be
timelike or spacelike, there is essentially only one type of
non-null vector in neutral signature. Note that a null vector for $g
\in [g]$ is null for all $g \in [g]$, so nullness of a vector with
respect to a conformal structure makes sense.

\subsection{Non-null case}

Given a neutral four dimensional ASD conformal structure $(\M,[g])$
with a non-null conformal Killing
vector $K$, the three dimensional space $\W$ of trajectories of $K$
inherits a conformal
structure $[h]$ of signature $(++-)$, due to
(\ref{KV}). The ASD condition on $[g]$ results in extra geometrical
structure on $(\W,h)$; it becomes a Lorentzian
\emph{Einstein-Weyl} space. This is called the Jones-Tod
construction, and is described in Section \ref{jonestod}. The next section
is an summary of Einstein-Weyl geometry.

\subsubsection{Einstein-Weyl geometry}
Let ${\cal W}$ be a three--dimensional manifold.\\
Given a conformal structure $[h]$ of signature\footnote{The
formalism in this section works in general dimension and signature
but we specialize to the case we encounter later.} $(2,1)$, a connection $D$ is said to
preserve $[h]$ if
\be \label{weylderiv}
Dh=\omega \otimes h,
\ee
for some $h \in [h]$, and a 1-form $\omega$. It is clear that
if (\ref{weylderiv}) holds
for a single $h \in [h]$ it holds for all, where $\omega$ will depend
on the particular $h \in [h]$.
(\ref{weylderiv}) is a natural condition; it is the requirement that null
geodesics of any $h \in [h]$ are also geodesics of $D$.

Given $D$ we can define its Riemann and Ricci curvature tensors
$W^{i}_{\ jkl}$, $W_{ij}$ in
the usual way. The notion of a curvature scalar must be modified,
because there is no distinguished metric in the conformal class to
contract $W_{ij}$ with. Given some $h \in [h]$ we can form $W=h^{ij}W_{ij}$.
Under a conformal transformation $h \rightarrow \phi^{2} h$, $W$
transformes as $W \rightarrow \phi^{-2} W$. This is because $W_{ij}$
unaffected by any conformal rescaling, being formed entirely out of
the connection $D$. $W$
is an example of a conformally weighted function, with weight $-2$.

One can now define a conformally invariant analogue of the Einstein
equation as follows:
\be \label{ew}
W_{(ij)}-\frac{1}{3} W h_{ij} = 0.
\ee
This are the \emph{Einstein-Weyl equations}\index{Einstein-Weyl equations}. 
Notice that the left hand side is well
defined tensor (i.e. weight $0$), since the weights of $W$ and $h_{ij}$
cancel. Equation (\ref{ew}) is the Einstein-Weyl equation for $(D,[h])$. It says
that given any $h \in [h]$, the Ricci tensor of $W$ is tracefree when
one defines the trace using $h$. Notice also that $W_{ij}$ is not
necessarily symmetric, unlike the Ricci-tensor for a Levi-Civita
connection.

In the special case that $D$ is the Levi-Civita connection of some
metric $h \in [h]$, (\ref{ew}) reduces to the Einstein equation. This
happens when $\omega$ is exact, because under $h \rightarrow
\phi^{2}h$, we get $\omega \rightarrow \omega + 2 d (\text{ln} \phi)$,
so if $\omega$ is exact a suitable choice of $\phi$ will transform
it to $0$, giving $D h = 0$ in (\ref{weylderiv}). All Einstein metrics
in 2+1 or 3 dimensions are spaces of constant curvature. The Einstein--Weyl
condition allows non-trivial degrees of freedom. The general solution
to (\ref{ew}) depends on four arbitrary functions of two variables.

In what follows, we refer to an Einstein-Weyl structure by
$(h,\omega)$. The connection $D$ is fully determined by this data
using $(\ref{weylderiv})$.

\subsubsection{Reduction by a non-null Killing vector; the Jones-Tod construction} \label{jonestod}
The Jones-Tod construction\index{Jones-Tod construction} relates
ASD conformal structures in four
dimensions to Einstein-Weyl structures in three
dimensions. In neutral signature it can be formulated as follows:
\begin{theorem}{\cite{jonestod}} \label{JT}
    Let ($\M,[g]$) be a
    neutral  ASD four manifold
    with a non-null conformal Killing vector $K$. An Einstein-Weyl
    structure on the space $\W$ of trajectories of $K$ is defined by
    \be \label{reducedEW}
    h:=|K|^{-2} g - |K|^{-4} \bK \odot \bK, \ \  \omega = 2
    |K|^{-2}\ast_{g}(\bK \wedge d \bK),
    \ee
    where $|K|^{2}:=g(K,K)$, $\bK:=g(K,.)$, and $\ast_{g}$ is the
    Hodge-$\ast$ of $g$. All EW structures arise in this
    way. Conversely, let $(h,\omega)$ be a three dimensional
    Lorentzian EW
    structure on $\W$, and let $(V,\eta)$ be a function and a
    1-form on $\W$ satisfying the generalised monopole 
equation\index{generalised monopole equation}
    \be
\label{general_mon}
    \ast_{h}(dV + \frac{1}{2} \omega V) = d \eta,
    \ee
    where $\ast_{h}$ is the Hodge-$\ast$ of $h$. Then
    $$
    g = V^{2} h - (d\phi + \eta)^{2}
    $$
    is a neutral ASD metric with non-null Killing vector $\p_{\phi}$.
\end{theorem}
This is a local theorem, so we may assume $\W$ is a manifold. A
vector in $\W$ is a vector field Lie-derived along the corresponding
trajectory in $\M$, and one applies the formulae (\ref{reducedEW}) to
this vector field to obtain $([h],\omega)$ on $\W$. In the
Riemannian case it has been successfully applied globally in certain
nice cases \cite{lebrun}. When one performs a conformal
transformation of $g$, one obtains a conformal transformation of $h$
and the required transformation of $\omega$, so this is a theorem
about conformal structures, though we have phrased it in terms of
particular metrics.

The Jones--Tod construction
was originally discovered using twistor theory in
\cite{jonestod}; since then other purely differential-geometric proofs have
appeared \cite{joyce,CP00}; although these are in
Riemannian signature the arguments carry over to the neutral case. In
Section \ref{twistor} we explain the twistorial argument that originally
motivated the theorem.

\vspace{10pt}

\subsubsection{Integrable systems and the Calderbank--Pedersen construction}
Applying the Jones--Tod correspondence to the special ASD conditions
discussed in Section 2 will yield special integrable systems in 2+1
dimensions. In each case of interest we shall assume that the symmetry
preserves the special geometric structure in four dimensions. This will give
rise to special Einstein--Weyl backgrounds, together with general solutions
of the generalised monopole equation (\ref{general_mon})
on these backgrounds. We can then seek
special monopoles such that the resulting ASD structure is conformal
to pseudo--hyper--K\"ahler.

An elegant framework for this is provided by the Calderbank--Pedersen
construction \cite{CP00}. In this construction
self--dual  complex (or null) structures on ${\cal M}$ correspond to
to shear-free  geodesic congruences (SFGC) on ${\cal W}$. This gives rise
to a classification of three-dimensional EW spaces according to
the properties of associated congruences. Below we shall list
the resulting reductions and integrable systems. In each case
we shall specify the properties of the associated congruence without
going into the details of the Calderbank--Pedersen correspondence.

\vspace{10pt}
\textbf{Scalar--flat K\"ahler with symmetry.
The $SU(\infty)$-Toda equation}. \ \
Let $({\cal M}, g)$ be a scalar--flat K\"ahler metric in neutral signature,
with a  symmetry  $K$  Lie deriving the K\"ahler form $\omega$.
 One can follow the steps
of LeBrun \cite{lebrun} to reduce the problem to a pair of coupled PDEs:
the  $SU(\infty)$-Toda equation\index{$SU(\infty)$-Toda equation} and 
its linearisation. The key step in the
construction is to use the moment map for $K$ as one of the coordinates, i .e.
define a function $t:\M\longrightarrow \R$ by $dt=K\hook \omega$. Then
$x, y$ arise as isothermal coordinates on two dimensional surfaces orthogonal
to $K$ and $d t$. The metric takes the form

\be \label{todametric}
g= V(e^u(dx^2+d y^2)-dt^2)-\frac{1}{V}(d\phi +\eta)^2,
\ee
where the function $u$ satisfies the $SU(\infty)$-Toda
equation
\be \label{toda}
(e^{u})_{tt} - u_{xx} - u_{yy} = 0,
\ee
and $V$ is a solution to its linearization -- the generalised monopole
equation   (\ref{general_mon}) .
The corresponding EW space from the Jones-Tod construction is
\be \label{todaEW}
h = e^{u}(dx^{2} +dy^{2})-dt^2, \ \ \omega = 2u_t dt.
\ee
It was shown in \cite{tod1} that the EW spaces that can be put in the
form (\ref{toda}) are precisely those possessing a shear-free
twist-free geodesic congruence. Given  the Toda EW space, any
solution to the monopole equation will yield a $(+ + - -)$
scalar  flat K\"ahler metric.
The special solution $V=cu_t$, where $c$ is a constant, will lead to a pseudo
hyperK\"ahler metric with symmetry. 

In \cite{nutku} 
solutions to  (\ref{toda}) were used to construct neutral 
ASD Ricci flat metrics without symmetries. 
\vspace{10pt}

\textbf{ASD Null K\"ahler with symmetry. The dKP equation}. \ \
Let $(\M, g, N)$ be an ASD null Kahler structure with a Killing vector
$K$ such that ${\cal L}_K N=0$.  In \cite{nullkahler} it was demonstrated that
there exist smooth real valued functions
$H=H(x, y, t)$ and $W=W(x, y, t)$  such that
\be
\label{scalarflat}
g=W_x(d y^2-4 d xd t-4H_xd t^2)-W_x^{-1}(d \phi-W_xd y-2W_yd t)^2
\ee
is an  ASD null K\"ahler metric on a circle bundle
${\cal M}\rightarrow {\cal W}$ if
\be
\label{Heqn}
H_{yy}-H_{xt} +H_xH_{xx}=0,
\ee
\be
\label{lindKP}
W_{yy}-W_{xt}+(H_xW_x)_x=0.
\ee
All real analytic
ASD null  K\"ahler metrics with symmetry
arise from this construction.

With definition $u=H_x$ the  $x$ derivative of equation
(\ref{Heqn}) becomes
$$(u_t-uu_x)_x=u_{yy},$$ which is
the dispersionless 
Kadomtsev--Petviashvili 
equation\index{dKP equation} (dKP)  originally
used in \cite{DMT}. The corresponding Einstein--Weyl structure is
$$
h = dy^{2} - 4 dx dt - 4 u dt^{2}, \ \ \omega=-4u_{x} dt.
$$
This EW structure possesses a covariant
constant null vector with weight $-\frac{1}{2}$, and in fact every
such EW structure with this property can be put into the above form. The
covariant constancy is with respect to a derivative on weighted
vectors that preserves their weight. Details can be found in \cite{DMT}.

The linear equation (\ref{lindKP}) is a (derivative of) the generalised
monopole equation from the Jones--Tod construction. Given a dKP Einstein--Weyl
structure, any solution to this monopole equation will yield and ASD Null
Kahler structure in four dimensions. The special monopole
$V= H_x/2$ will yield a pseudo--hyper--K\"ahler structure with symmetry whose
self--dual derivative is null.

\vspace{10pt}

\textbf{Pseudo-hypercomplex with symmetry. 
The hyper--CR equation\index{hyper--CR equation}.} \ \
Let us assume that a pseudo--hyper-complex four manifold admints a symmetry
which Lie derives all (pseudo) complex structures.  This implies \cite{duna:hydro}
that the EW structure is locally given by
\[
h=(d y+ud t)^2-4( d x+wd t)d t,\qquad
\omega =u_xd y+(uu_x+2u_y)d t,
\]
where $u(x, y, t)$ and $w(x, y, t)$ satisfy a system  of quasi-linear PDEs
\be
\label{PMA}
u_t+w_y+uw_x-wu_x=0,\qquad u_y+w_x=0.
\ee
The corresponding pseudo--hypecomplex metric will arise form any solution
to this coupled system, and its linearisation (the generalised
monopole (\ref{general_mon})). The special monopole $V=u_x/2$ leads to
pseudo-hyper-K\"ahler metric with triholomorphic homothety.

\vspace{10pt}

\textbf{Interpolating integrable system.} \ \
There exists a  dispersionless integrable system
which contains the dKP 
equation
and the hyper--CR equation as  special cases.
This {\em interpolating integrable system} is given by
\cite{interpol}
\be
\label{uw_eq}
u_y+w_x=0, \qquad u_t+w_y-c(uw_x-wu_x)+buu_x=0,
\ee
where $b$ and $c$ are constants and $u, w$ are smooth functions of
$(x,y,t)$.  It admits a Lax pair 
\begin{eqnarray*}
L_0&=&\frac{\p}{\p t}+(cw+bu-\lambda cu-\lambda^2)\frac{\p}{\p x}+
b(w_x-\lambda u_x)\frac{\p}{\p \lambda}\\ 
L_1&=&\frac{\p}{\p y}-(cu +\lambda)\frac{\p}{\p x}-bu_x\frac{\p}{\p\lambda}.
\end{eqnarray*}

Setting $b=0, c=-1$ gives the hyper--CR equation  
and setting $c=0, b=1$  gives the dKP equation. 
In fact one constant can always be eliminated from (\ref{uw_eq}) 
by redefining the coordinates and it is only the ratio of $b/c$ which 
remains. The interpolating integrable system arises as the general 
symmetry reduction of pseudo--hyperk\"ahler metric  with a conformal Killing vector whose self--dual derivative is null \cite{interpol}.

\vspace{10pt}

\subsection{Null case} \label{null}

Given a neutral four dimensional ASD conformal structure $(\M,[g])$
with a null conformal Killing
vector\index{null conformal Killing
vector} $K$, the three dimensional space of trajectories of $K$
inherits a degenerate conformal structure of signature $(+-0)$, and
the Jones-Tod construction does not hold. The situation was
investigated in detail in \cite{DW} and \cite{cald06}.
It was shown that $K$ defines a pair of
totally null foliations of $\M$, one by $\alpha$-surface and one by
$\beta$-surfaces; these foliations intersect along integral curves of
$K$ which are null geodesics.
In spinors, if $K^a = \iota^A o^{A'}$ then an $\alpha$-plane
distribution is defined by $o^{A'}$, and a $\beta$-plane
distribution by $\iota^A$, and it follows from the Killing equation
that these distributions are integrable.

The main result from \cite{DW} is that there is a canonically defined
\emph{projective structure}\index{projective structure} on 
the two-dimensional space of
$\beta$-surfaces $U$ which arises as a quotient of ${\cal M}$ by a
distribution $\iota^A\e_{AA'}$. A more general framework
where the distribution $\iota^A\e_{AA'}$ is still integrable, but
$\iota^A o^{A'}$ is not a symmetry for any $o^{A'}\in\Gamma (S')$
was recently developed by
Calderbank \cite{cald06} and extended by Nakata \cite{nakata2}.

A projective structure is an
equivalence class of connections, where two connections are
equivalent if they have the same unparameterized geodesics. In
Section \ref{twistor} we will explain the twistor theory that led to
the observation that projective structures are involved, and give a
new example of a twistor construction.

It turns out that one can explicitly write down all ASD conformal
structures with null conformal Killing vectors in terms of their underlying
projective structures as follows:

\begin{theorem}\cite{DW} \label{localform}
Let $(\M,[g],K)$ be a smooth neutral signature ASD conformal structure
with null conformal Killing vector. Then there exist
local coordinates $(\phi,x,y,z)$ and $g \in [g]$ such that $K=\p_\phi$ and
$g$ has one
of the following two forms, according to whether the  twist
$\bK \wedge d \bK$
vanishes or not $(\bK := g(K,.))$:
\begin{enumerate}
\item $\bK \wedge d \bK =0$.
\begin{multline} \label{nontwistinggeneral}
g = (d\phi + (z A_3 - Q) dy)(dy - \beta dx) - \\(dz - (z (-\beta_y +
A_1 + \beta A_2 + \beta^2 A_3))dx - (z(A_2 + 2 \beta A_3)+ P)dy)dx,
\end{multline}
where $A_1,A_2,A_3,\beta,Q,P$ are arbitrary functions of $(x,y)$.

\item $\bK \wedge d \bK \neq 0$.
\begin{multline} \label{twistinggeneral}
g = (d\phi + A_3\p_zG  dy + (A_2\p_z G  + 2 A_3 (z \p_zG - G)
- \p_z\p_yG)dx)(dy - z dx) \\  - \p_z^2G dx(dz - (A_0 + z A_1 + z^2
A_2 + z^3 A_3)dx),
\end{multline}
where $A_0,A_1,A_2,A_3$ are arbitrary functions of $(x,y)$, and
$G$ is a function of $(x,y,z)$ satisfying the following PDE:
\begin{equation} \label{Gequation}
(\p_x + z \p_y + (A_0+zA_1+z^2 A_2 + z^3 A_3) \p_z) \p_z ^2G = 0.
\end{equation}
\end{enumerate}
\end{theorem}
The functions $A_i(x,y)$ in the metrics (\ref{nontwistinggeneral}) and
(\ref{twistinggeneral}) determine projective structures on the
two dimensional space $U$ in the following way. A two projective structure
in two dimensions is equivalent to a second-order ODE
\begin{equation} \label{secondorderode}
\frac{d^2 y}{dx^2} = A_3(x,y) \Big( \frac{d y}{dx} \Big)^3 +
A_2(x,y) \Big( \frac{dy}{dx} \Big)^2 + A_1(x,y) \Big( \frac{dy}{dx}
\Big) + A_0(x,y),
\end{equation}
obtained by choosing local coordinates $(x,y)$ and eliminating the
affine parameter from the geodesic equation. The $A_{i}$ functions can be
expressed in terms of combinations of connection coefficients that
are invariant under projective transformation.
In (\ref{twistinggeneral}) all the $A_i, i=0,1,2,3$ functions occur
explicitly in the metric. In (\ref{nontwistinggeneral}) the function
$A_0$ does not explicitly occur. It is determined by the following
equation:
\be
\label{beta_equation}
A_0 = \beta_x + \beta \beta_y - \beta A_1 - \beta^2 A_2 - \beta^3 A_3.
\ee
If the projective structure is flat, i.e.$A_i=0$ and $\beta=P=0$ then
(\ref{nontwistinggeneral})
is Ricci flat \cite{Pl75}, and in fact this is the most general
ASD Ricci flat metric with a null Killing vector which preserves
the pseudo hyperK\"ahler structure \cite{BGRPP}.
More generally, if the
projective structure comes from a Riemannian metric on $U$ then there 
will always exist
a (pseudo-)K\"ahler structure in the conformal class (\ref{twistinggeneral})
if $G=z^2/2+\gamma(x, y)z+\delta(x, y)$
for certain  $\gamma, \delta$  \cite{BDE08,DT09}.

It is interesting that integrable systems are not involved in the
null case, given their ubiquity in the non-null case.

\section{Twistor theory} \label{twistor}
In Riemannian signature, given an ASD conformal structure $(\M,[g])$ in four
dimensions one can form a 2-sphere bundle over it, and endow this with
an integrable complex structure by virtue of anti-self-duality \cite{AHS}.
The resulting complex manifold $\PT$ is called the 
twistor space\index{twistor space}. The
original manifold is the moduli space of rational curves in $\PT$
preserved under a certain anti-holomorphic involution\index{anti-holomorphic involution}, and 
one can
recover the conformal structure by looking at how the rational
curves intersect one another. Hence the $(\M,[g])$ is completely
encoded in $\PT$ and its anti-holomorphic involution. The important
feature of a successful twistor construction is that the original geometry
becomes encoded in the holomorphic geometry of the twistor space, and
can be recovered from this.

Neutral signature ASD conformal structures cannot be encoded purely
in holomorphic geometry as in the Riemannian case. This is not
surprising as generically they are
not analytic. However, there is a recent twistor construction due to
LeBrun-Mason \cite{lebrunmason} in the
neutral case that uses a mixture of holomorphic and smooth
ingredients; we review this in Section \ref{lebrunmason}. Let us now
review the differences in Riemannian and neutral signature.

In the Riemannian case, if one expresses the metric in terms of a
null tetrad as in (\ref{themetric}) then the basis vectors $\e_{AA'}$
must be complex, as there are no real null vectors. The spin bundles
are complex two dimensional vector bundles $\mathbb{S}$,
$\mathbb{S'}$, with an isomorphism $T_{\mathbb{C}} \M \cong
\mathbb{S} \otimes \mathbb{S'}$, at least locally. One then takes the
projective bundle $\mathbb{P} \mathbb{S'}$, which has $\CP$ fibres. Even
if $\mathbb{S'}$ does not exist globally, the bundle $\mathbb{P} \mathbb{S'}$
does exist globally, since the $\mathbb{Z}_2$ obstruction to existence of
a spin bundle is eliminated on projectivizing. Concretely, $\mathbb{PS}'$
is the bundle of complex self-dual totally null 2-planes; from this
description it clearly exists globally.

Now one can form the
$L_{A}$ vectors as in Theorem \ref{penrose}, where now $\pi^{A'}$ are
complex (the homogeneous fibre coordinates of $\mathbb{P} \mathbb{S'}$). The
connection coefficients in the expression for $\et_{AA'}$ will now be
complex, and satisfy certain Hermiticity properties that we need not
go into. The $L_{A}$ span a two complex dimensional distribution
on the complexified tangent space of $\mathbb{P} \mathbb{S'}$, and the
Riemannian version of Theorem \ref{penrose} is that this 
distribution is complex
integrable iff the metric is ASD. Together with $\p_{\bar{\lambda}}$,
where $\lambda$ is the inhomogeneous fibre coordinate on $P
\mathbb{S'}$, we obtain a complex three dimensional distribution
$\Pi$, satisfying $\Pi \cap \bar{\Pi} = 0$. If the metric is
ASD, $\Pi$ is complex integrable and defines a complex structure
on $\mathbb{P} \mathbb{S'}$. This construction works globally. It
was discovered by Atiyah, Hitchin and Singer \cite{AHS}.

In the neutral case one can complexify the real spin bundles $S$,
$S'$ and obtain $T_{\mathbb{C}} \M \cong S_{\C}¥ \otimes
S'_{\C}$ as in the Riemannian case. One can define a complex
distribution $\Pi$
distribution on $P S'_{\C}$, by allowing $\pi^{A'}$ in Theorem
\ref{penrose} to be complex. The key point is that the vectors
$L_{A}$ become totally real when $\pi^{A'}$ is real. So on the
hypersurface $P S' \subset P S'_{\C}$, the distribution $\text{span} \{ \Pi,
\p_{\bar{\lambda}} \}$ no longer satisfies $\Pi \cap \bar{\Pi}=0$, so
does not define an almost complex structure. When $\pi^{A'}$ is not real,
the distribution spanned by
$L_{A}$ and $\p_{\bar{\lambda}}$ does define an almost
complex structure, which is integrable when $g$ is ASD. We obtain
two non-compact regions in $P S'_{\C}$, each of which possesses an integrable
complex structure, separated by a hypersurface $PS'$. This is more  complicated
than the Riemannian case, where the
end result is simply a complex manifold.
Nevertheless, the construction is reversible in a precise sense given
by  Theorem  of LeBrun-Mason which we review in Section
\ref{lebrunmason} (Theorem \ref{machine}).

Before describing the work of LeBrun-Mason we review the analytic
case, where one can complexify and work in the holomorphic category.

\subsection{The analytic case} \label{analytic}
In this section we work locally. Standard references for this
material are \cite{wardwells,hitchin}.

Suppose a neutral four dimensional ASD conformal
structure $(\M,[g])$
 is analytic in some coordinate system. Then we can
complexify by letting the coordinates become complex variables, and
we obtain a holomorphic conformal structure $(\M_{\C},[g_{\C}])$. If
each coordinate is defined in some connected open set on $\R$, then
one thickens this slightly on both sides of the axis to obtain a
region in $\C$ on which the complex coordinate is defined. The
holomorphic conformal structure is obtained by picking a real metric
$g$ and allowing the coordinates to be complex to obtain $g_{\C}$. Then
$[g_{\C}]$ is the equivalent class of $g_{\C}$ up to multiplication
by nonzero holomorphic functions.

>From Theorem \ref{penrose}, which is valid equally for holomorphic
metrics, we deduce that given any holomorphic
$\alpha$-plane at a point, there is a holomorphic $\alpha$-surface
through that point. Assuming we are working in a suitably
convex neighbourhood so that the space of such $\alpha$-surfaces is
Hausdorff, we define $\PT$ to be this space. $\PT$ is a three dimensional
complex manifold, since the space of $\alpha$-planes at a point is
one complex dimensional and each surface is complex codimension two
in $\M_{\C}$. This is summarised in  the double fibration
\index{double fibration} picture
\be
\label{doublefib}
{\cal M_{\C}}\stackrel{p}\longleftarrow
{ P S'_{\C}}\stackrel{q}\longrightarrow {\cal PT},
\ee
where $q$ is the quotient by the 
twistor distribution $L_A$\index{twistor distribution}.

If we had started with a Riemannian metric this would
lead to the same twistor space, locally, as the Atiyah-Hitchin-Singer construction
described above, though we shall not demonstrate this
here. A point $x \in \M$, corresponds to an embedded $\CP \subset \PT$,
since there is a $\CP$ of \as s through $x$. By varying the point $x
\in \M$ we obtain a four complex parameter family of $\CP$s.

$\PT$ inherits an anti-holomorphic involution $\sigma$. To
describe $\sigma$, note that there is an anti-holomorphic involution
$\tau$ of $\M_{\C}$ that fixes real points, i.e. points of $\M
\subset \M_{\C}$. This is just the map from a coordinate to its
complex conjugate, so we can arrange our complexification
regions in which the coordinates are defined so that $\tau$ maps the
regions to themselves. Now $\tau$ will map holomorphic \as s to
holomorphic \as s, so gives an anti-holomorphic involution $\sigma$
on $\PT$. One way to see this is to note that \as s are totally
geodesic as the geodesic shear free condition
\[
\pi^{A'}\pi^{B'}\nabla_{AA'}\pi_{B'}=0
\]
is equivalent to $C_{A'B'C'D'}$,
and consider the holomorphic geodesic
equation. Using the fact that the connection coefficients are real, one can show
that the involution $\tau$ will map the null geodesics in an \as \ to
other null geodesics in another \as. The \as s fixed by this are the
real \as s in $\M$.

In terms of $\PT$, this last fact means that $\sigma$ fixes an equator of each of
the four complex parameter family of embedded $\CP$s. Moreover, an
\as \ through a real point gets mapped to one through that same point
since the point is fixed by $\tau$. So
the $\CP$s that are fixed by $\sigma$ are a four real parameter
family corresponding to $\M$, we call these \emph{real} $\CP$ s.

How does one recover the neutral conformal structure from the data
$(\PT,\sigma)$? As described above, $\M$ is the moduli space of $\CP$s
fixed by $\sigma$. Now a vector at a point in $\M$ corresponds to a
holomorphic section\index{holomorphic sections} of the normal bundle $\OO(1) \oplus \OO(1)$ of
the corresponding real $\CP$
in $\PT$, such that the section `points' to another real $\CP$. We
define a vector to be null if this holomorphic section has a zero.
Since vanishing of a section of
$\OO(1) \oplus \OO(1)$ is a quadratic condition, this gives a
conformal structure. One can prove that this conformal structure is
ASD  fairly easily, by showing that the required \as s must exist in
terms of the holomorphic geometry.

Moreover, special conditions on a $g_{\C} \in [g_{\C}]$ can be encoded
into the holomorphic geometry of the twistor space:
\begin{itemize}
\item Holomorphic fibration $\theta:\PT\rightarrow \CP$
  $\longleftrightarrow$
corresponds to
hyper-hermitian conformal structures \cite{boyer,duna:twisted}.
\item
Preferred section of $\kappa^{-1/2}$  which
vanishes at exactly two points on each twistor line corresponds to
scalar--flat K\"ahler $g_{\C}$ \cite{pontecorvo}.
\item
Preferred section of $\kappa^{-1/4}$ corresponds
to ASD null K\"ahler
$g_{\C}$ \cite{nullkahler}.
\item
Holomorphic fibration $\theta:\PT\rightarrow \CP$ and holomorphic
isomorphism \[\theta^*{\cal O}(-4)
\cong\kappa\] correspond to hyper-K\"ahler
$g_{\C}$
\cite{penrose,AHS,hitchin}.
\end{itemize}
Here $\kappa$ is a holomorphic canonical bundle\index{canonical bundle} of $\PT$, and ${\cal O}(-4)$
is a power of the tautological bundle on the base of $\theta$.
To obtain a real metrics the structures above must be preserved by
an anti-holomorphic  involutions
fixing a real equator of  each rational curve in $\PT$.

%

It is worth saying a few words about the construction of solutions of
integrable systems using the twistor correspondence. It is shown in
Section \ref{symmetries} that a number of well-known integrable
systems  $2+1$  dimensions are special
cases of ASD conformal structures. Analytic solutions to these
integrable systems therefore
correspond\footnote{This correspondence is not
one-one due to coordinate freedom.}
to twistor spaces $\PT$.
There will be extra conditions on
$\PT$, depending on the special case in question. However, solutions
to the integrable systems are not always analytic.

\subsubsection{Symmetries and twistor spaces} \label{symtwistor}
In Section \ref{symmetries} we discussed the appearance of
Einstein-Weyl structures and projective structures in the cases of a
non-null and null Killing vector respectively. In both cases twistor
theory was the key factor in revealing these correspondences. We
shall now explain this briefly. In \cite{hitchin}, Hitchin gave three twistor
correspondences. He considered complex manifolds containing embedded
$\CP$s with normal bundles $\OO(1)$, $\OO(2)$ and $\OO(1) \oplus
\OO(1)$ respectively. Kodaira 
deformation theory\index{Kodaira deformation theory} guarantees a local
moduli space of embedded $\CP$s, whose complex dimension is the
dimension of the space of holomorphic sections of the corresponding
normal bundle, i.e. 2, 3, 4 respectively. By examining how nearby
curves intersect, he deduced that the moduli space inherits a
holomorphic projective structure, Einstein-Weyl structure, or ASD conformal
structure respectively. He also showed that the construction is
reversible in each case.

Now given a four dimensional holomorphic ASD conformal structure, its
twistor space is the space of $\alpha$-surfaces, as described in
Section \ref{analytic}. A conformal Killing vector preserves the
conformal structure, so preserves $\alpha$-surfaces, giving a
holomorphic vector field on the twistor space.

If the Killing vector
is non-null then the vector field on twistor space $\PT$ is non-vanishing.
This is because the Killing vector is transverse to any
$\alpha$-surface, as it is non-null. In this case one can quotient
the three dimensional twistor space by the induced vector field, and
it can be shown \cite{jonestod} that the resulting two dimensional complex
manifold contains $\CP$s with normal bundle $\OO(2)$. Using Hitchin's
results, this corresponds to a three dimensional Einstein-Weyl structure. This the
twistorial version of the Jones-Tod construction, Theorem \ref{JT}.

If the Killing vector is null then the induced vector field on the
twistor space $\PT$ vanishes on a hypersurface. This is because at each
point, the Killing vector is tangent to a single $\alpha$ surface.
Hence it preserves a foliation by $\alpha$-surfaces, and vanishes at
the hypersurface in twistor space corresponding to this foliation.
However, one can show \cite{DW} that it is possible to continue the vector
field on twistor space to a one-dimensional distribution ${\hat{K}}$
that is
nowhere vanishing. Quotienting $\PT$ by this distribution gives a
two-dimensional complex manifold $Z$ containing $\CP$s with normal bundle
$\OO(1)$. Using Hitchin's results, this corresponds to a two
dimensional projective structure. This is the twistorial version of
the correspondence described in Section \ref{null}. The situation is
illustrated by the following diagram:

\begin{figure}
\label{Figg}
\centering
\includegraphics[height=8cm]{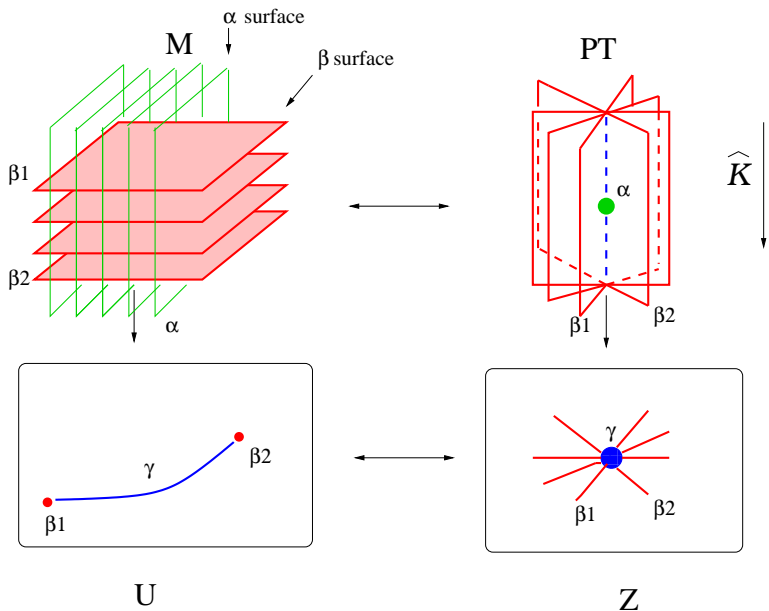}
\caption{Relationship between $\M$, $U$, $\PT$ and $Z$.}
\end{figure}

In $\M$, a one parameter family of
\bs \ is shown, each of which intersects a one parameter family of \as
s, also shown. The \bs s correspond to a projective structure geodesic
in $U$, shown at the bottom left.

The \bs s in $\M$ correspond to surfaces in $\PT$, as discussed
above. These surfaces intersect at the dotted line, which corresponds
to the one parameter family of \as s in $\M$. When we quotient $\PT$ by
$\hat{K}$ to get $Z$, the surfaces become twistor lines in $Z$, and the
dotted line becomes a point at which the twistor lines intersect; this
is shown on the bottom right. This family of twistor lines
intersecting at a point corresponds to the geodesic of the projective
structure.

\vspace{5pt}
\textbf{Example}\footnote{We thank Paul Tod for his help with this example.}.
Here we give an explicit construction of the twistor space of an
analytic neutral ASD conformal structure with a null Killing vector,
from the reduced projective structure twistor space. We take $Z$ to
be the total space of $\OO(1)$. This is the twistor space of the flat
projective structure. Now suppose we are given a 1-form
$\omega$ on $U$. We shall complexify the setup and regard $\omega$ as
holomorphic a holomorphic connection on a holomorphic line
bundle $B \rightarrow U$. This gives rise to a holomorphic line bundle
$E \rightarrow Z$, where the
vector space over $z \in Z$ is the space of parallel sections of $B$
over the geodesic in $U$ corresponding to $z$. The \emph{twistor
lines} in $Z$ are the two-parameter family of embedded $\CP$s, each
corresponding to the set of geodesics through a single point in $U$.
We denote the twistor line corresponding to a point $x \in U$ by
$\hat{x}$. Now $E$ restricted to a twistor line $\hat{x}$ is trivial,
because to specify a parallel section of $B$ through any geodesic
through $x$, one need only know its value at $x$. This is a simple
analogue of the Ward correspondence relating solutions of the
anti-self-dual Yang-Mills equations on $\C^{4}$ to vector bundles
over the total space of $\OO(1) \oplus \OO(1)$ that are trivial on
twistor lines. The situation here is simpler since there are no
PDEs involved; this is because there are no integrability conditions
for a space of parallel sections to exist on a line. As with the Ward
correspondence, the construction is reversible, i.e. given a
holomorphic line bundle trivial on twistor lines one can find a
connection on $U$ to which it corresponds in the manner described
above. We will not prove this here, it is simply a case of mimicking
the argument for the Ward correspondence \cite{wardwells}.

Now to create the
twistor space $\PT$, we must tensor $E$ with a line bundle $L$ so
that $E \otimes L$ restricts to $\OO(1)$ on the twistor lines in $Z$.
Then the total space of $E \otimes L$ will have embedded $\CP$s with
normal bundle $\OO(1) \oplus \OO(1)$, so will be a twistor space for
an ASD conformal structure. For $L$ we choose the pull back of
$\OO(1)$ to the total space of $\OO(1)$.

Let us now make the above explicit. Let $\lam$, $\tl$ be the
inhomogeneous coordinate on the two patches $U_{0}$, $U_{1}$ of
$\CP$. The total space of $\OO(1)$ can be coordinatized as follows.
Let $\mu$ be the fibre coordinate over $U_{0}$, and $\tm$ the
fibre coordinate over $U_{1}$. The line bundle transition relation on
the overlap is $\tm = \frac{1}{\lambda} \mu$.

Now suppose we have a line bundle $E \rightarrow Z=\OO(1)$, that is
trivial on holomorphic sections of $Z\rightarrow \CP$.
Let $\tau$, $\ttau$ be the fibre
coordinates on the two patches, satisfying a transition relation
$\ttau = F(\lam,\mu) \tau$, where $F(\lam,\mu)$ is holomorphic and
nonvanishing on the overlap, i.e. for $\lam \in \C - \{ 0 \}$, $\mu
\in \C$. In sheaf terms, $F$ is an element of
$H^{1}(\OO(1),\OO^{*})$. Now the short exact sequence
\be \label{SES}
0 \rightarrow \Z \rightarrow \OO \rightarrow \OO^{*} \rightarrow 0
\ee
gives rise to a long exact sequence, part of which is:
\be \label{LES}
\ldots \rightarrow H^{1}(\OO(1),\Z) \rightarrow
H^{1}(\OO(1),\OO) \rightarrow H^{1}(\OO(1),\OO^{*}) \rightarrow
H^{2}(\OO(1),\Z)\rightarrow \ldots
\ee
The first term in (\ref{LES}) vanishes and the final term
is $\Z$, by topological considerations. The final term gives the
Chern class of the line bundle determined by the element of
$H^{1}(\OO(1),\OO^{*})$. This vanishes for $E$, since it is trivial
on twistor lines. The third arrow in (\ref{SES}) is the exponential
map. Together these facts
imply that $F$ can be written $F(\lambda,\mu) =
e^{f(\lambda,\mu)}$, where $f(\lambda,\mu)$ is a holomorphic function
on the overlap that may have zeros. After twisting by $L$, we obtain
the following transition function for $E \otimes L$, again using
$\tau$, $\ttau$ as fibre coordinates:
\be \label{transition}
\ttau =  \frac{1}{\lambda} e^{f(\lambda, \mu)} \tau.
\ee
To find the conformal structure we must find the four parameter
family of twistor lines in $E\otimes L$.
The two parameter family in $\OO(1)$ is
given in one patch by $\mu(\lambda) = X \lambda + Y$, and in the other
by
$\tm(\tl) = X + \tl Y$. Restricting to one of
these we can split $f$:
\be
\label{fsplit}
f(\lambda, X \lambda + Y) = h(X,Y,\lambda) - \tilde{h}(X,Y,1/ \lambda),
\ee
where $h$ and $\tilde{h}$ are functions on $U\times \CP$ holomorphic in
$\lambda$ and $1 / \lambda$ respectively.
For fixed $(X,Y)$ there is then a further two
parameter family of twistor lines, given by
\be \label{tau0}
    \tau(\lam) = e^{-h(X,Y,\lambda)}(W -\lam Z)
\ee
in one patch, and
\be \label{tau1}
\ttau (\tl) = e^{-\tilde{h}(X,Y,\tl)}(\tl W - Z).
\ee
It is easy to check that (\ref{transition}) is satisfied by
(\ref{tau0}) and (\ref{tau1}).

One must now calculate the conformal structure on the moduli space of lines
parametrised by $X^a=(X, Y, W, Z)$  by determining the
quadratic condition for a section of the normal bundle to a twistor
line to vanish. The sections of the normal bundle to $\hat{x}\subset \PT$
correspond to tangent vectors in $T_x{\cal M}$, and sections with one zero
will determine null vectors and therefore the conformal structure.

Using the identity $(\p_X-\lambda\p_Y)f=0$
together with (\ref{fsplit}) we deduce  (by Liouville theorem or using power series) that
\be
\label{gexpression}
\lambda\frac{\p h}{\p Y} - \frac{\p h}{\p X} = \lambda B(X,Y)
- A(X,Y),
\ee
for some analytic functions $A$,$B$.

Now take the variation of $\mu(\lambda)$ and $\tau(\lambda)$
for a small change  $\delta X^a$ to obtain
\begin{eqnarray}
0 = \delta \mu  &=& \dY + \lambda \dX,  \label{muvariation} \\
0 = \delta \tau &=& e^{-h}(-\frac{\p h }{\p X} \dX - \frac{\p h}{\p Y}
\dY)(W -\lambda Z) + e^{-h}(\dW -\lambda  \dZ). \label{tauvariation}
\end{eqnarray}
Substituting  $\lambda=-\delta Y/\delta X$ from the first expression
to the second, using  (\ref{gexpression}) and multiplying the
resulting expression by $\delta X$
we find that the conformal structure is represented by the
following metric:
\be
g = dX dW + dY dZ - (W dX + Z dY)(A(X,Y) dX + B(X,Y) dY).
\ee
This conformal structure possesses the null conformal Killing vector $K = W
\p_{W} + Z \p_{Z}$, which is twisting. The global holomorphic
vector field on $\PT$ induced by $K$ is
$\tau \p_{\tau} = \ttau \p_{\ttau}$ where
the equality holds on the intersection of the two coordinate patches.
This vanishes on
the hypersurface
defined by $\tau = 0$ in one patch and $\ttau = 0$ in the other,
which intersects each twistor line at a single point,
as we expect from the argument
in Section \ref{symtwistor}. The 1-form $\omega=AdX+BdY$ in $g$ is the inverse
Ward transform of $F\in H^1({\cal O}(1), {\cal O}^*)$.

To compare with
(\ref{twistinggeneral}) one must transform to coordinates $(\phi,x,y,z)$
in which $K=\p_{\phi}$. Dividing by a conformal factor $W$,
transforming  with $(\phi,x,y,z)=(\text{log}
W,Y,-X,Z/W)$, and then translating $\phi$ to eliminate an arbitrary
one function of $(x, y)$ gives
\be
\label{finalconfstructure}
g=(d\phi + f(x,y) dx)(dy - z dx) - dz dx,
\ee
a special case of (\ref{twistinggeneral}) with flat projective
structure, and $G=z^2/2-zC(x, y)$, where $f=\p_y C$.

If we take the coordinates to be
real we obtain a neutral metric. The twistor space $\PT$ fibres
over $Z=\OO(1)$ and this fibres over $\CP$,
so $\PT$ fibers over $\CP$ and
(\ref{finalconfstructure}) is pseudo-hyperhermitian.

To construct an example of a
conformal structure with \emph{non-twisting} null
Killing vector one uses an affine line bundle over $Z=\OO(1)$; see
\cite{DW} for details.

\subsection{LeBrun--Mason construction} \label{lebrunmason}
Here we describe recent work of LeBrun and
Mason\index{LeBrun--Mason construction} in which a general, global twistor
construction is given for neutral metrics. We will only be able to give a
crude paraphrase, and refer the reader to the original paper \cite{lebrunmason} for
details. Note that their paper uses the opposite duality conventions to
ours; they use self-dual conformal structures with integrable $\beta$-plane
distributions.

We described above how a neutral ASD conformal structure $(\M,[g])$
gives rise to a complex structure on $\CP$ bundle over $\M$, which
degenerates on a hypersurface. The following theorem of LeBrun-Mason
is a converse to this, and is the closest one can come to a general
twistor construction in the neutral case:

\begin{theorem}{\cite{lebrunmason}}\label{machine}
 Let $\M$ be a smooth connected
$4$-manifold,  and let $\varpi: {\mathcal X}\to \M$ be a smooth
$\CP$-bundle. Let   $\varrho  : {\mathcal X}\to {\mathcal X}$  be an involution
which
commutes with   $\varpi$, and has  as fixed-point  set
${\mathcal X}_\varrho$  an $S^1$-bundle over $\M$ which disconnects
${\mathcal X}$ into two closed $2$-disk bundles ${\mathcal X}_\pm$
with common boundary ${\mathcal X}_\varrho$.
Suppose that $\Pi \subset T_\C {\mathcal X}$ is
a    distribution of complex $3$-planes on ${\mathcal X}$
such that
\begin{enumerate}
\item
$\varrho _* \Pi = \overline{\Pi}$;
\item   the restriction of $\Pi$ to ${\mathcal X}_+$
is smooth
and  involutive,
\item $\Pi \cap \overline{\Pi} = 0$ on   ${\mathcal X}-{\mathcal X}_\varrho$,
\item $\Pi \cap \ker \varpi_*$  is the $(0,1)$ tangent space of the $\CP$ fibers
of $\varpi$,
\item  the restriction of $\Pi$ to a fiber of ${\mathcal X}$  has
first Chern class $-4$ with respect to the complex orientation.
\end{enumerate}
Then ${\cal E}=\Pi\cap T{\mathcal X}_\varrho$ is an integrable distribution of real   $2$-planes on
${\mathcal X}_\varrho$, and  $\M$ admits a
  unique smooth split-signature  ASD conformal structure $[g]$
  for which the $\alpha$-surfaces are the projections
  via $\varpi$  of the integral manifolds of ${\cal E}$.
  \end{theorem}
This theorem provides a global twistor construction for neutral ASD
four manifolds, whereas the analytic construction of the last section
only works locally.

At first sight the theorem does not seem like a promising method of
generating ASD conformal structures, since the conditions required on
the $\CP$ bundle over $\M$ are complicated, and it is not clear how
one might construct examples. This obstacle is overcome in \cite{lebrunmason} by
deforming a simple example
(another example was given by Nakata \cite{Nakata}).

Consider the conformally flat
neutral metric $g_0$ given by (\ref{standard})
on $\M = S^{2} \times S^{2}$ that is just the
difference of the standard sphere metrics on each factor. The underlying
manifold
$\M$ can be realised as the space of $\CP$s embedded in
$\mathbb{CP}^{3}$ that are invariant under the complex conjugate
involution, which we call the \emph{real} $\CP$s; these real $\CP$s
are the fibres of the bundle
$\mathcal{X}$. The involution $\varrho$ of $\mathcal{X}$ is induced
by the complex conjugate involution of $\mathbb{CP}^{3}$. The fixed
point set $\mathcal{X}_{\varrho}$
consists of the invariant equators of the real $\CP$s, and is
therefore a circle bundle over $\M$. The closed disc bundles
$\mathcal{X}_{\pm}$ are obtained by slicing the real $\CP$s at their
invariant equator, and throwing away one of the open halves. The
fixed point set of the complex conjugate involution is the standard
embedding of $\mathbb{RP}^{3}$, and this is the space of
\as s in $\M = S^{2} \times S^{2}$. Taking all the real $\CP$s
through a point $p \in \mathbb{RP}^{3}$ gives an \as.

To obtain the $\Pi$ from  Theorem \ref{machine},
take $\mathcal{X}_{+}$ and construct a map $f$ to $\mathbb{CP}^{3}$ as follows. On
the interior of $\mathcal{X}_{+}$, $f$ is a diffeomorphism onto
$\mathbb{CP}^{3}-\mathbb{RP}^{3}$. The boundary $\p
\mathcal{X}_{+}$ gets mapped by $f$ to
$\mathbb{RP}^{3}$, by taking a point in
$\p \mathcal{X}_{+}$, i.e. a holomorphic disc and a point $p$ on
$\mathbb{RP}^{3}$ lying on the intersection of the disc with
$\mathbb{RP}^{3}$, to the point $p$. Let
$f_{*}^{1,0}: T_{\C} \mathcal{X}_{+} \rightarrow T^{1,0}
\mathbb{CP}^{3}$ be the $(1,0)$ part of the derivative of $f$. Then
the $\Pi$ of  Theorem  \ref{machine} is defined on $\mathcal{X}_{+}$ by
$$
\Pi = \text{ker} f_{*}^{1,0} \subset T_{\C} \mathcal{X}_{+}.
$$
Note that $f$ maps the five dimensional boundary $\p \mathcal{X}_{+}$
to the three dimensional space $\mathbb{RP}^{3}$; this means that on
the boundary $\Pi$ restricts to the complexification of a real two-plane
distribution, direct summed with the complexification of the
direction into the disc.
The $\Pi$ here agrees with the one described at the
beginning of Section \ref{twistor}, defined in terms of the twistor
distribution $L_{A}$ on $P S'_{\C}$.

The idea for creating new spaces satisfying the conditions of Theorem
\ref{machine}
is to deform the standard embedding of $\mathbb{RP}^{3}
\subset \mathbb{CP}^{3}$ slightly. One forms the bundle
$\mathcal{X}_{+}$ by taking the same $S^{2} \times S^{2}$ family of
real holomorphic discs, now with
boundary on the deformed embedding of $\mathbb{RP}^{3}$. The complex
distribution $\Pi$ is formed in the same way as in the conformally
flat case described above. One then patches two copies of
$\mathcal{X}_{+}$ together to form the bundle $\mathcal{X}$, and it
is shown that this satisfies the conditions of the Theorem. It is
also shown that the resulting
conformal structure on $S^{2} \times S^{2}$ has the property that all
null geodesics are embedded circles; conformal structures with this
property are termed \emph{Zollfrei}. It turns out that all ASD
conformal structures close enough in a suitable sense to the
conformally flat one are Zollfrei, and in fact the twistor description
gives a complete understanding of ASD conformal structures near the
standard one. The embedded $\mathbb{RP}^{3}$ is the real twistor
space, i.e. the space of \as s in $\M$, and a significant portion of
\cite{lebrunmason} is devoted to showing that the for a space-time oriented
Zollfrei 4-manifold the real twistor space must be $\mathbb{RP}^{3}$,
making contact with the picture of a deformed $\mathbb{RP}^{3}
\subset \mathbb{CP}^{3}$.

We mention that there is another twistor--like construction of
smooth ASD conformal structures with avoids the holomorphic methods 
altogether \cite{grossman}. In this approach one views the real twistor
curves in $\mathbb{RP}^3$ as 
solutions to a system of two second order nonlinear
ODEs. The ODEs have to satisfy certain conditions (expressed in terms
of point invariants) if their solution spaces are equipped 
with ASD conformal structures.

\section{Global results} \label{global}
In the last section we outlined the global twistor construction for
neutral ASD four manifolds due to LeBrun-Mason, which they used to
construct Zollfrei metrics on $S^{2} \times S^{2}$. In this section
we review the known explicit constructions of globally defined
neutral ASD conformal structures on various compact and non-compact
manifolds.

\subsection{Topological restrictions}

Existence of a neutral metric on a four manifold $\M$ imposes
topological restrictions on $\M$. A neutral inner product on a
four dimensional vector $V$ space splits $V$ into a direct sum $V =
V_+ \oplus V_-$, where the inner product is positive definite
on $V_+$ and negative definite on $V_-$. So a neutral metric $g$ on a four
manifold $\M$ splits the tangent bundle
\be \label{tangentsplitting}
T \M = T_+ \M \oplus T_-\M,
\ee
where $T_{\pm}$ are two dimensional subbundles of $T\M$.
Conversely, given such a splitting one can construct neutral
metrics on $\M$ by taking a difference of positive definite metrics on
the vector bundles $T_+ \M$ and $T_- \M$.

If $\M$ admits a non-vanishing 2-plane field ${\cal E}$
(a real two--dimensional
distribution), then a splitting of the
form (\ref{tangentsplitting}) can be found by taking ${\cal E}$
to define $T_+ \M$,
choosing a Riemannian metric, and letting $T_- \M$ be the orthogonal
complement. So a four manifold $\M$ admits a neutral metric iff it
admits a 2-plane field.

The topological conditions for existence of a 2-plane field were discovered
by Hirzebruch and Hopf and are as follows:
\begin{theorem}{\cite{HH}}
A compact smooth four--manifold
$\M$ admits a field of 2-planes iff $\tau[\M]$ and $\chi[\M]$ satisfy a pair of
conditions
\begin{eqnarray*}
3 \tau[\M] + 2 \chi [\M] \in \Omega(\M), \\
3 \tau[\M] - 2 \chi [\M] \in \Omega(\M),
\end{eqnarray*}
where
$$
\Omega[\M] = \{ \mu_{\M}(w,w) \in \mathbb{Z}: \text{$w$ are arbitrary
elements in $H^2(\M,\mathbb{Z})/\text{Tor}$}  \} .
$$
\end{theorem}
Here $\tau[\M]$, $\chi[\M]$ are the signature and Euler characteristic
respectively, and $\mu_{\M}$ is the intersection form on
$H^2(\M,\mathbb{Z})/\text{Tor}$.

A neutral metric implies that the structure group of the tangent bundle
can be reduced to $O(2,2)$, by choosing orthonormal bases in each patch.
In fact $O(2,2)$ has four connected components, so there are various
different orientability requirements one can impose. The simplest is to
require the structure group to reduce to the identity component
$SO_+(2,2)$. It is shown in \cite{ML} that this is equivalent to the
existence of a field of \emph{oriented} 2-planes, i.e. an orientable two
dimensional sub-bundle of the tangent bundle. The topological restrictions
imposed by this were discovered by Atiyah:
\begin{theorem}\cite{atiyah}
Let $\M$ be a compact oriented smooth manifold of dimension $4$, such that
there exists a field of \emph{oriented} 2-planes on $\M$. Then
\be \label{atiyahconditions}
\chi[\M] \equiv 0 \ \text{mod} \ 2, \ \ \ \ \ \chi[\M] \equiv \tau[\M] \
\text{mod} \ 4.
\ee
\end{theorem}
In fact Matsushita showed \cite{mat:SC} that for a simply-connected
4-manifold, (\ref{atiyahconditions}) are actually sufficient for the
existence of an oriented field of 2-planes.
A more subtle problem is to determine
topological obstructions arising from existence of an ASD neutral
metric. This deserves further study.

\subsection{Tod's scalar-flat K\"ahler metrics on $S^{2} \times
  S^{2}$}
\label{Tod_sub}
Consider $S^{2} \times S^{2}$ with the conformally flat metric
described in Sections  \ref{conf_comp} and  \ref{lebrunmason},
i.e. the difference of the
standard sphere metrics on each factor. Thinking of each sphere as
$\CP$ and letting $\zeta$ and $\chi$ be non-homogeneous coordinates for
the spheres, this metric is given by (\ref{standard}).
As we have already said,  $g_0$ is scalar flat,
indefinite K\"ahler. The obvious complex structure $J$ gives a closed two
form  and $\Omega:=g_0(J.,.)$. Moreover $g_0$  clearly
has a high degree of symmetry,
since the 2-sphere metrics have rotational symmetry. In \cite{tod},
Tod found deformations of $g_0$ preserving the scalar-flat K\"ahler
property, by using the explicit expression (\ref{todametric}) for
neutral scalar-flat K\"ahler metrics with symmetry. Take the explicit
solution
$$
e^{u} = 4 \frac{1-t^{2}}{(1+x^{2}+y^{2})^{2}}
$$
to (\ref{toda}), which can be obtained by demanding $u =
f_1(x,y) + f_2(t)$. There remains a linear equation for $V$. Setting
$W=V(1-t^{2})$  and  performing
the coordinate transformation $t = \cos \theta$, $\zeta = x+ i y$
gives
\be \label{sfksym}
g = 4 W \frac{d\zeta d \bar{\zeta}}{(1+\zeta \bar{\zeta})^{2}} - W d
\theta^{2} - \frac{\sin^{2}\theta}{W}(d \phi + \eta)^{2},
\ee
and $W$ must solve a linear equation. This metric reduces to
(\ref{standard}) for $W=1$, $\eta = 0$, with $\theta$, $\phi$
standard coordinates for the second sphere. Tod shows that on
differentiating the linear equation for $W$ and setting $Q = \frac{\p
W}{\p t}$, one obtains the ultrahyperbolic wave equation
\be \label{uhwave}
\nabla_{1}^{2} Q = \nabla_{2}^{2}Q,
\ee
where $\nabla_{1,2}$ are the Laplacians on the 2-spheres, and $Q$ is
independent of $\phi$, i.e. is axisymmetric for one of the sphere
angles. Equation (\ref{uhwave}) can be
solved using Legendre polynomials, and one obtains non-conformally
flat deformations of (\ref{standard}) in this way. In the process one
must check that $W$ behaves in such a way that (\ref{sfksym}) extends over
$S^{2} \times S^{2}$.

The problem of relating these explicit metrics to the Zollfrei metrics
on $S^2 \times S^2$ known to exist by results described in Section
\ref{lebrunmason} appears to be open.

In a recent paper, Kamada \cite{kam:sfk} rediscovered the above metrics,
and showed that a compact neutral scalar-flat K\"ahler manifold with a
Hamiltonian $S^1$ symmetry must in fact be $S^{2} \times S^{2}$. Here a
Hamiltonian $S^1$ symmetry is an $S^1$ action preserving the K\"ahler form,
and which possesses a moment map. In the case of $S^{2} \times S^{2}$ case,
there is always a moment map since the manifold is simply connected. Without
the symmetry, there are other neutral scalar-flat K\"ahler manifolds. For
example, take a Riemann surface $\Sigma$ with a constant curvature metric $g$.
Then on $\Sigma \times \Sigma$, the metric $\rho_1^* g - \rho_2^* g$ is neutral
scalar-flat K\"ahler, where $\rho_i$ are the projections onto the first and
second factors.

\subsection{Compact neutral hyperk\"ahler metrics}
The only compact four
dimensional Riemannian hyperk\"ahler manifolds are the complex
torus with the flat metric and $K3$ with a Ricci-flat 
Calabi-Yau\index{Calabi-Yau}
metric. In the neutral case, Kamada showed in \cite{kam:kodaira} that
a compact pseudo-hyperk\"ahler four manifold must be either a complex
torus or a primary Kodaira surface. In the complex torus case, the
metric need not be flat, in contrast to the Riemannian case.
Moreover in both cases one can write down explicit non-flat examples, in
contrast to the Riemannian case where no explicit non-flat Calabi-Yau
metric is known.

To write down explicit examples, consider the following hyperk\"ahler
metric
\be \label{ppwave}
g = d\phi dy - dz dx - Q(x,y) dy^{2},
\ee
for $Q$ and arbitrary function. This is the neutral version of the
pp-wave metric of general
relativity \cite{Pl75},
and is a special case of (\ref{nontwistinggeneral}), where the
underlying projective structure is flat. It is
non-conformally flat for generic $Q$. Define complex coordinates
$z_{1} = \phi + i z$, $z_{2} = x + i y$ on $\C^{2}$. By quotienting the
$z_{1}$ and $z_{2}$ planes by lattices one obtains a product of
elliptic curves, a special type of complex torus.
If we require $Q$ to be periodic with respect to the $z_{2}$ lattice,
then (\ref{ppwave}) descends to a metric on this manifold. Likewise,
a primary Kodaira surface can be
obtained as a quotient of $\C^{2}$ by a subgroup of the group of
affine transformations, and again by assuming suitable periodicity in
$Q$ the metric (\ref{ppwave}) descends to the quotient.
In our framework \cite{DW}
we compactify the flat projective space $\R^2$ to two--dimensional torus
$U=T^2$ with the projective structure coming from the flat metric. Both
$\phi$ and $z$ in (\ref{ppwave})
 are taken to be periodic, thus leading
to $\hat{\pi}:{\cal M}\longrightarrow U$, the holomorphic
toric fibration\index{toric fibration} over a torus. Assume the suitable
periodicity on the function $Q:U\longrightarrow \R$.
This leads to a commutative diagram
\begin{eqnarray*}
&\M&\\
T^2&\downarrow&\searrow \hat{\pi}^*Q\\
&U& \stackrel{Q}{\longrightarrow}\R .
\end{eqnarray*}
In the framework of  \cite{kam:kodaira} and \cite{Fino}
the K\"ahler structure on ${\cal M}$ is given by
$\omega_{flat}+i\p\overline{\p}(\hat{\pi}^*Q)$, where $(\p, \omega_{flat})$
is the flat K\"ahler structure on the Kodaira surface induced from $\C^2$.

As remarked in \cite{kam:kodaira}, the existence of
pseudo-hyperk\"ahler metrics on complex tori other than a product of
elliptic curves is an open problem.

\subsection{Ooguri-Vafa metrics\index{Ooguri-Vafa metrics}}
In \cite{oogurivafa} Ooguri, Vafa and Yau constructed a class of
non-compact neutral hyper-K\"ahler metrics
on cotangent bundles of Riemann surfaces with genus $\geq 1$,
using the Heavenly equation formalism. This is similar to (\ref{heavy1}),
but one takes a different $(+ + - -)$
real section of ${\cal M}_\C$. Instead of using the real
coordinates we set
\[
w=\zeta, \quad y=\bar{\zeta}, \quad z=ip, \quad x=-i\bar{p}, \qquad \zeta,p\in \C
\]
with $\Omega=i(p\bar{\zeta}-\bar{p}\zeta)$ corresponding to the flat metric.
Let $\Sigma$ be a Riemann surface with a local  holomorphic
coordinare $\zeta$,
such that the K\"ahler metric on $\Sigma$ is
$h_{\zeta\bar{\zeta}}d\zeta d\bar{\zeta}$.  Suppose that
$p$ is a  local complex coordinate for fibres of the cotangent
bundle $T^* \Sigma$.  If $\omega$ is the K\"ahler form for
a neutral metric $g$ then $g_{i\bar{j}}=\p_i \p_{\bar{j}}\Omega$
for a function
$\Omega$ on the cotangent bundle. Then the equation
$$\text{det} \ g_{i\bar{j}} = -1$$
is equivalent to the first Heavenly equation (\ref{heavy1}), and gives
a Ricci-flat ASD neutral metric.

The idea in \cite{oogurivafa} is to suppose that $\Omega$ depends only on the
globally defined function $X = h^{\zeta \bar{\zeta}} p \bar{p}$, which
is  the length of the cotangent vector
corresponding to $p$.
There is a globally defined holomorphic $(2, 0)$-form
$\omega = d\zeta \wedge dp$,
which is the holomorphic part of the standard symplectic form
on the cotangent bundle, so  $(\zeta, p)$ are the holomorpic
coordinates in the Pleba\'nski
coordinate system. The heavenly equation reduces to
an ODE for $\Omega(X)$ and Ooguri-Vafa show that
for solutions of this ODE to exist $h$ must have
constant negative curvature, so $\Sigma$ has genus greater than
one. In this case  one can solve the ODE to find
\[
\Omega=2\sqrt{A^2+BX
}+A\ln{\frac{\sqrt{A^2+BX}-A}{\sqrt{A^2+BX}+A}},
\]
where $A, B$ are arbitrary positive constants.
The metric $g$ is well behaved when $X\rightarrow 0$ (or $p\rightarrow 0$),
as in this limit
$\Omega\rightarrow \ln{(X)}$ and $g$  restricts to
$h_{\zeta\bar{\zeta}}d\zeta d\bar{\zeta}$ on $\Sigma$ and
$-h^{\zeta\bar{\zeta}}dp d\bar{p}$ on the fibres. In
the limit $X\rightarrow\infty$ the
metric is flat. To see it one needs to chose a uniformising coordinate
$\tau$ on $\Sigma$ so that $h$ is a metric on the upper half plane.
Then make a coordinate transformation $\zeta_1=\tau\sqrt{p},
\zeta_2=\sqrt{p}$. The holomorphic two form is still $d \zeta_1\wedge d
\zeta_2$, and the K\"ahler potential
$\Omega=i(\zeta_2\bar{\zeta}_1-\zeta_1\bar{\zeta}_2)\sqrt{B}$ yields
the flat metric.

Ooguri-Vafa  also observed that the pp-wave metric (\ref{ppwave})
can be put onto
$T^*\Sigma$, by requiring $Q(x,y)$ to satisfy certain symmetries.
Globally defined neutral metrics on non-compact manifolds were also
studied by Kamada and Machida in \cite{KM}, where they obtained many neutral
analogues of well known ASD Riemannian metrics.


\begin{thebibliography}{10}

\bibitem{ashtekar}
{\sc Ashtekar, A., Jacobson, T., and Smolin, L.}
\newblock A new characterization of half-flat solutions to einstein's equation.
\newblock {\em Comm. Math. Phys. 115\/} (1988), 631--648.

\bibitem{atiyah}
{\sc Atiyah, M.~F.}
\newblock Vector fields on manifolds.
\newblock {\em Arbeitsgemeinschaft f\"ur Forschung des Landes Nordrhein
  Westfalen 200\/} (1970), 7--24.

\bibitem{AHS}
{\sc Atiyah, M.~F., Hitchin, N.~J., and Singer, I.~M.}
\newblock Self-duality in four-dimensional {R}iemannian geometry.
\newblock {\em Proc. Roy. Soc. London Ser. A 362\/} (1978), 425--461.

\bibitem{BGRPP}
{\sc Barrett, J., Gibbons, G.~W., Perry, M., and Pope, C.}
\newblock Kleinian geometry and the $n=2$ superstring.
\newblock {\em Int. J. Mod. Phys. A 9\/} (1994), 1457--1494.

\bibitem{boyer}
{\sc Boyer, C.~P.}
\newblock A note on hyperhermitian four-manifolds.
\newblock {\em Proc. Amer. Math. Soc. 102\/} (1988), 157--164.

\bibitem{osserman}
{\sc Brozoz-Vazquez, M., Garcia-Rio, E., and Vazquez-Lorenzo, R.}
\newblock Conformally osserman four-manifolds whose conformal jacobi operators
  have complex eigenvalues.
\newblock {\em Proc. Roy. Soc. London Ser. A 462\/} (2006), 1425--1441.


\bibitem{druma}
{\sc Bogdanov, L. V., Dryuma, V.S., and Manakov S.V.}
\newblock
On the dressing method for {D}unajski anti-self-duality equation.
\newblock{\em nlin/0612046\/}.

\bibitem{B00}
{\sc Bryant, R.~L.}
\newblock Pseudo-{R}iemannian metrics with parallel spinor fields and vanishing
  {R}icci tensor.
\newblock In {\em Global analysis and harmonic analysis\/} (2000), vol.~4 of
  {\em Semin. Congr.}, Soc. Math. France, Paris, pp.~53--94.

\bibitem{BDE08} {\sc Bryant, R. L., Dunajski, M. and Eastwood, M.}
\newblock Metrisability of two-dimensional projective structures.
\newblock  {\em J. Differential Geometry. 83\/} (2009), 465-499.


\bibitem{cald06}
{\sc Calderbank, D. M.~J.}
\newblock Selfdual 4-manifolds, projective structures, and the
  {D}unajski-{W}est construction.
\newblock {\em math.DG/0606754\/}.

\bibitem{CP00}
{\sc Calderbank, D. M.~J., and Pedersen, H.}
\newblock Self-dual spaces with complex structures, {E}instein-{W}eyl geometry
  and geodesics.
\newblock {\em Ann. Inst. Fourier 50\/} (2000), 921--963.

\bibitem{DJS}
{\sc Dancer, A.~S., Jorgensen, H.~R., and Swann, A.~F.}
\newblock Metric geometries over the split quaternions.
\newblock {\em Rend. Sem. Mat. Univ. Politec. Torino 63\/} (2005), 119--139.

\bibitem{gueo}
{\sc Davidov, J.,  Grantcharov, G.,  Mushkarov, O. and  Yotov, M.}
\newblock  Para--hyperhermitian surfaces. 
\newblock {\em Bull. Math. Soc. Sci. Math. Roumanie (N.S.) 52(100)\/} 
(2009) 281--289. 

\bibitem{Derdzinski}
{\sc Derdzi\'nski, A.}
\newblock Self-dual {K}\"ahler manifolds and {E}instein manifolds of dimension
  four.
\newblock {\em Compositio Math. 49\/} (1983), 405--433.

\bibitem{duna:twisted}
{\sc Dunajski, M.}
\newblock The twisted photon associated to hyper-hermitian four-manifolds.
\newblock {\em J. Geom. Phys. 30\/} (1999), 266--281.

\bibitem{nullkahler}
{\sc Dunajski, M.}
\newblock Anti-self-dual four-manifolds with a parallel real spinor.
\newblock {\em Proc. Roy. Soc. London Ser. A 458\/} (2002), 1205--1222.

\bibitem{duna:hydro}
{\sc Dunajski, M.}
\newblock A class of {E}instein-{W}eyl spaces associated to an integrable
  system of hydrodynamic type.
\newblock {\em J. Geom. Phys. 51\/} (2004), 126--137.

\bibitem{DMT}
{\sc Dunajski, M., Mason, L.~J., and Tod, K.~P.}
\newblock Einstein-{W}eyl geometry, the d{K}{P} equation and twistor theory.
\newblock {\em J. Geom. Phys. 37\/} (2001), 63--93.

\bibitem{DW}
{\sc Dunajski, M., and West, S.}
\newblock Anti-self-dual conformal structures from projective structures.
\newblock {\em Comm. Math. Phys. 272\/} (2007), 85--118.

\bibitem{interpol} {\sc Dunajski, M}
\newblock An interpolating dispersionless integrable system.
\newblock {\em J. Phys. A: Math. Theor. 41\/} (2008) 315202.  


\bibitem{D_book} {\sc Dunajski, M.}
\newblock {\em Solitons, Instantons \& Twistors}
Oxford Graduate Texts in Mathematics vol. ~19, 
\newblock Oxford University Press, 2009.

\bibitem{DT09} {\sc Dunajski, M., Tod, K. P.} 
\newblock Four Dimensional Metrics Conformal to K\"ahler. 
\newblock {\em Math. Proc. Camb. Soc. 148\/} (2010) 485-503.



\bibitem{Fino}
{\sc Fino, A., Pedersen, H., Poon, Y.-S., and Sorensen, M.~W.}
\newblock Neutral calabi-yau structures on kodaira manifolds.
\newblock {\em Comm. Math. Phys. 248\/} (2004), 255--268.

\bibitem{grossman}
{\sc Grossman, D.}
\newblock Torsion-free path geometries and integrable 
second order ODE systems, 
\newblock{\em Sel. Math., New Ser. 6\/} 
(2000), 399-442.

\bibitem{HH}
{\sc Hirzebruch, F., and Hopf, H.}
\newblock Felder von {F}l\"achenelementen in 4-dimensionalen
  {M}annigfaltigkeiten.
\newblock {\em Math. Ann. 136\/} (1958), 156--172.

\bibitem{hitchin}
{\sc Hitchin, N.~J.}
\newblock Complex manifolds and {E}instein's equations.
\newblock In {\em Twistor geometry and Non-Linear Systems\/} (1982), Lecture
  Notes in Math., 970, Springer, pp.~73--99.

\bibitem{hitchin:HS}
{\sc Hitchin, N.~J.}
\newblock Hypersymplectic quotients.
\newblock {\em Acta Acad. Sci. Tauriensis 124\/} (1990), 169--180.

\bibitem{ivanov}
{\sc Ivanov, S., and Zamkovy, S.}
\newblock Parahermitian and paraquaternionic manifolds.
\newblock {\em Diff. Geom. App. 23\/} (2005), 205--234.

\bibitem{john}
{\sc John, F.}
\newblock The ultrahyperbolic wave equation with four independent variables.
\newblock {\em Duke Math. J. 4\/} (1938), 300--322.

\bibitem{jonestod}
{\sc Jones, P.~E., and Tod, K.~P.}
\newblock Minitwistor spaces and einstein-weyl spaces.
\newblock {\em Class. Quant. Grav. 2\/} (1985), 565--577.

\bibitem{joyce}
{\sc Joyce, D.~D.}
\newblock Explicit construction of self-dual $4$-manifolds.
\newblock {\em Duke Math. J. 77\/} (1995), 519--552.

\bibitem{kam:kodaira}
{\sc Kamada, H.}
\newblock Neutral hyperk\"ahler structures on primary {K}odaira surfaces.
\newblock {\em Tsukuba J. Math. 23\/} (1999), 321--332.

\bibitem{kam:sfk}
{\sc Kamada, H.}
\newblock Compact scalar-flat indefinite {K}\"ahler surfaces with {H}amiltonian
  ${S}^1$-symmetry.
\newblock {\em Comm. Math. Phys. 254\/} (2005), 23--44.

\bibitem{KM}
{\sc Kamada, H., and Machida, Y.}
\newblock Self-duality of metrics of type $(2,2)$ on four-dimensional
  manifolds.
\newblock {\em T\^ohoku Math. J. 49\/} (1997), 259--275.

\bibitem{kob}
{\sc Kobayashi, S., and Nomizu, K.}
\newblock {\em Foundations of differential geometry, {V}olume II}.
\newblock No.~15 in Interscience Tracts in Pure and Applied Mathematics. John
  Wiley \& Sons, 1969.

\bibitem{Law}
{\sc Law, P.}
\newblock Classification of the weyl curvature spinors of neutral metrics in
  four dimensions.
\newblock {\em J. Geom. Phys.\/} (2006).

\bibitem{lebrun}
{\sc LeBrun, C.}
\newblock Explicit self-dual metrics on $\mathbb{CP}^2 \# \ldots \#
  \mathbb{CP}^2$.
\newblock {\em J. Diff. Geom. 34\/} (1991), 223--253.

\bibitem{lebrunmason}
{\sc LeBrun, C., and Mason, L.~J.}
\newblock Nonlinear gravitons, null geodesics, and holomorphic disks.
\newblock {\em Duke Math. J. (To appear)\/}.


\bibitem{nutku}
{\sc  Malykh, A. A., Nutku, Y., and Sheftel, M. B.}
\newblock 
Lift of noninvariant solutions of heavenly equations 
from three to four dimensions and new ultra-hyperbolic metrics.
\newblock {\em arXiv:0704.3335v2\/}.

\bibitem{masonnewman}
{\sc Mason, L.~J., and Newman, E.~T.}
\newblock A connection between the einstein and yang-mills equations.
\newblock {\em Comm. Math. Phys. 121\/} (1989), 659--668.

\bibitem{masonwoodhouse}
{\sc Mason, L.~J., and Woodhouse, N. M.~J.}
\newblock {\em Integrability, self-duality, and twistor theory}, vol.~15 of
  {\em London Mathematical Society Monographs}.
\newblock Oxford University Press, 1996.

\bibitem{mat:SC}
{\sc Matsushita, Y.}
\newblock Fields of 2-planes on compact simply-connected smooth 4-manifolds.
\newblock {\em Math. Ann. 280\/} (1988), 687--689.

\bibitem{ML}
{\sc Matsushita, Y., and Law, P.}
\newblock Hitchin-{T}horpe type inequalities for pseudo-{R}iemannian
  4-manifolds of metric signature $(++--)$.
\newblock {\em Geom. Dedicata 87\/} (2001), 65--89.

\bibitem{MTbook}
{\sc Mimura, Y., and Takeno, H.}
\newblock {\em Wave geometry}.
\newblock Sci. Rep. Res. Inst. Theoret. Phys. Hiroshima Univ. No. $2$. 1962.

\bibitem{Nakata}
{\sc Nakata, F.}
\newblock Singular self--dual zollfrei metrics and 
twistor correspondence.
\newblock {\em math.DG/0607276\/}.

\bibitem{nakata2}
{\sc Nakata, F.}
\newblock
Self-dual Zollfrei conformal structures with alpha-surface foliation.
\newblock{\em math/0701116\/}.

\bibitem{oogurivafa}
{\sc Ooguri, H., and Vafa, C.}
\newblock Geometry of {N}$=2$ strings.
\newblock {\em Nucl. Phys. B 361\/} (1991), 469--518.

\bibitem{Park92}
{\sc Park, Q.~H.}
\newblock $2${D} sigma model approach to $4${D} instantons.
\newblock {\em Internat. J. Modern Phys. A 7\/} (1992), 1415--1447.

\bibitem{penrose}
{\sc Penrose, R.}
\newblock Nonlinear gravitons and curved twistor theory.
\newblock {\em Gen. Rel. Grav. 7\/} (1976), 31--52.

\bibitem{PR86}
{\sc Penrose, R., and Rindler, W.}
\newblock {\em Spinors and {S}pace-time 1 \& 2}.
\newblock Cambridge Monographs on Mathematical Physics. Cambridge University
  Press, 1988.

\bibitem{Pl75}
{\sc Pleba\'nski, J.~F.}
\newblock Some solutions of complex {E}instein equations.
\newblock {\em J. Math. Phys. 16\/} (1975), 2395--2402.

\bibitem{PR76}
{\sc Pleba\'nski, J.~F., and Robinson, I.}
\newblock Left-degenerate vacuum metrics.
\newblock {\em Phys. Rev. Lett. 37\/} (1976), 493--495.

\bibitem{pontecorvo}
{\sc Pontecorvo, M.}
\newblock On twistor spaces of anti-self-dual hermitian surfaces.
\newblock {\em Trans. Am. Math. Soc. 331\/} (1992), 653--661.

\bibitem{old_pleban}
{\sc Sibata, T., and Morinaga, K.}
\newblock A complete and simpler treatment of wave geom.
\newblock {\em Proc. Amer. Math. Soc. 102\/} (1988), 157--164.

\bibitem{taubes}
{\sc Taubes, C.~H.}
\newblock The existence of anti-self-dual conformal structures.
\newblock {\em J. Diff. Geom. 36\/} (1992), 163--253.

\bibitem{tod1}
{\sc Tod, K.~P.}
\newblock Scalar-flat {K}\"ahler and hyper-{K}\"ahler metrics from
  painl\'eve-iii.
\newblock {\em Class. Quant. Grav. 12\/} (1992), 1535--1547.

\bibitem{tod}
{\sc Tod, K.~P.}
\newblock Indefinite conformally-{ASD} metrics on ${S}^2 \times {S}^2$.
\newblock In {\em Further Advances in Twistor Theory\/} (2001), vol.~III,
  Chapman \& Hall/CRC, pp.~61--63.
\newblock Reprinted from \emph{Twistor Newsletter 36}, 1993.

\bibitem{Walker}
{\sc Walker, A.~G.}
\newblock Canonical form for a riemannian space with a parallel field of null
  planes.
\newblock {\em Quart. J. Math. Oxford. 1\/} (1950), 69--79.

\bibitem{wardwells}
{\sc Ward, R.~S., and Wells, R.~O.}
\newblock {\em Twistor geometry and field theory}.
\newblock Cambridge Monographs on Mathematical Physics. Cambridge University
  Press, 1990.

\bibitem{Wo92}
{\sc Woodhouse, N. M.~J.}
\newblock Contour integrals for the ultrahyperbolic wave equation.
\newblock {\em Proc. Roy. Soc. London Ser. A 438\/} (1992), 197--206.

\end{thebibliography}

\end{document}